%% file: Deformation_Classification_of_Quartic_Surfaces_with_Simple_Singularities.tex
\documentclass{amsart}

\usepackage{hyperref}
\usepackage{amsmath,color,graphicx,amssymb,amscd}
\usepackage{degt}
\usepackage{tikz-cd} 
\usepackage{booktabs}
\usepackage{comment} 
\usepackage[paper=portrait,pagesize]{typearea}
\usepackage{longtable}
\usepackage{array, boldline, makecell, booktabs}

\newlength\mylength

\def\L{\mathbf{L}}
\def\K{\mathcal{K}}

\def\X{\mathcal{X}}

\def\Z{\mathbb{Z}}
\def\Q{\mathbb{Q}}

\def\disc{\operatorname{disc}}

\def\Aut{\operatorname{Aut}}

\def\ie{\emph{i.e\PERIOD}}
\def\eg{\emph{e.g\PERIOD}}
\def\iq.{\emph{i.q\PERIOD}}
\def\cf.{\emph{cf.}}
\def\viz.{\emph{viz\PERIOD}}
\def\vs.{\emph{vs\PERIOD}}

\begin{document}

\title[Classification of Quartic Surfaces with Simple Singularities]{Deformation Classification of Quartic Surfaces with Simple Singularities}

\author[\c{C}.~G\"{u}ne\c{s} ~Akta\c{s} ]{\c{C}\.{i}sem G\"{u}ne\c{s} Akta\c{s}}
\address{
 Department of Engineering Sciences, Abdullah G\"{u}l University\\
 38080, Kayseri, Turkey} \email{cisem.gunesaktas@agu.edu.tr}

\thanks{The author was supported by the \textit{T\"{U}B\.{I}TAK} grant $118$F$413$}



\date{}

\dedicatory{}

\begin{abstract}
We give a complete equisingular deformation classification of simple spatial quartic surfaces which are in fact $K3$-surfaces.
\end{abstract}

\maketitle


\section{Introduction}\label{introduction}
Throughout the paper, all varieties are over the field $\mathbb{C}$ of complex numbers.

\subsection{Motivation and historical remarks}

Thanks to the global Torelli theorem~\cite{K3} and surjectivity of the period map~\cite{periodmap}, the equisingular deformation classification of singular projective models of K3-surfaces with any given polarization becomes a mere computation. The most popular models studied intensively in the literature are plane sextic curves and spatial quartic surfaces. Using the arithmetical reduction~\cite{Alex1}, Akyol and Degtyarev~\cite{Alex2} completed the problem of equisingular deformation classification of simple plane sextics. Simple quartic surfaces which play the same role in the realm of spatial surfaces as sextics do for curves, are a relatively new subject, promising interesting discoveries. We confine ourselves to \textit{simple} quartics only, \textit{i.e.}, those with $\mathbf{A}$--$\mathbf{D}$--$\mathbf{E}$ type singularities
(quartic surfaces with a non-simple singular point, \ie, when the quartic is not a K3-surface, are quite different, see Degtyarev \cite{Alex89a}, \cite{Alex94b}). The work was originated in Urabe~\cite{Urabe3, Urabe2, Urabe1} and extended by Yang~\cite{Yang} who gave a complete list of sets of singularities realized by simple quartics. Then, after a period of oblivion, it was resumed by G\"{u}ne\c{s} Akta\c{s}~\cite{Cisem1} where she  obtained the classification of the so-called \emph{nonspecial} simple quartics by using the same aproach as in Degtyarev and Akyol~\cite{Alex2}.
In the meanwhile, Shimada~\cite{Shimada.connEllK3}, inspired by the work on quartics, produced a complete list of the connected components of the moduli space of Jacobian elliptic K3-surfaces (which can be regarded as $U$-polarized K3-surfaces).

It has become quite apparent that the classical approach to quartics and sextics based on the defining equations is bound to fail (see, \eg, Artal et al. \cite{Artal2001}, \cite{Artal2002}, \cite{Artal2003}; Degtyarev \cite{Alex94a}, or Oka and Pho \cite{OkaPo2002}); even when the classification is already known, it requires a tremendous amount of work to find the defining equations of curves/surfaces with large sets of singularities. On the other hand, the more modern $K3$-theoretic approach, pioneered by Urabe~\cite{Urabe3, Urabe2} and Yang~\cite{Yang}, has already demonstrated its productivity, see, \eg ~\cite{Alex2},~\cite{Cisem1} or~\cite{Shimada.connEllK3}. 

Also worth mentioning is the vast literature related to the study of the deformation classification problems in the \emph{real} case, singular or smooth: the classification of real smooth quartic surfaces by Kharlamov~\cite{Kharlamov}, the study on moduli space of real $K3$-surfaces  by Nikulin~\cite{Niku4} or the
work on quartic spectrahedra by Degtyarev and Itenberg~\cite{AI} and Ottem \emph{et al.}~\cite{Sturmfels}.
\subsection{Principal results}\label{principal.results}
The present paper originates from the paper \cite{Cisem1} where I have made a contribution to the systematic study of simple spatial quartics. Our principal result is extending the classification given in~\cite{Cisem1} for only nonspecial simple quartics to the whole space of simple quartics and, thus, completing the equisingular deformation classification of simple quartic surfaces. This result closes a long standing  project initiated by Persson~\cite{Persson} and Urabe~\cite{Urabe3, Urabe2, Urabe1}.

A set of simple singularities of a simple quartic can be identified with a root lattice (see \cite{Dufree} and \S\ref{configuration}) and recall that the total Milnor number $\mu(X)$ of a simple quartic $X\subset \mathbb{P}^3$ is given by the rank of the corresponding root  lattice. One has $\mu(X)\leq 19$ (see \cite{Urabe2}, \emph{cf.}, \cite{Persson}); if $\mu(X)=19$, the quartic $X$ is called  \textit{maximizing}. (Recall that maximizing quartics are projectively rigid.) 

Our classification is based on the arithmetical reduction found in Degtyarev and Itenberg~\cite{AI} reducing the equisingular deformation classification of simple quartic surfaces to a purely lattice-theoretical question. The resulting arithmetical problem is solved by first using Nikulin’s theory of discriminant forms~\cite{Niku2}. Then, the computation is done separately for the maximizing ($\mu(X)= 19$) and non-maximizing ($\mu(X)\leq 18$) case; for the former, we use Gauss’s theory of binary quadratic form~\cite{Gauss}; for the latter, apply Miranda-Morrison’s results~\cite{MM1,MM2,MM3}, reducing the analysis of the orthogonal groups of indefinite lattices to a relatively simple computation in finite abelian groups.

Our classification is obtained by implementing all the algorithms given in \S\ref{algorithm.for.classification} in GAP~\cite{GAP} as the number of classes is huge (about 12,000). The original code used in ~\cite{Cisem1} where we settled the case of nonspecial quartics was based partially on manually precomputed data specific to a particular degree. 
Therefore, the code has been extended to implement a complete version of the Miranda--Morrison \cite{MM1,MM2,MM3} and a genuinely repetition free enumeration of realizable configurations (see \S \ref{configuration}). The principal novelty, compared to the nonspecial case where $\K=0$ , is that the imprimitivity, where  $\K\ne0$, is also taken into account; furthermore,  we have developed a more advanced method of computing the Miranda–-Morrison homomorphism based on lifting reflections to a p-adic lattice (see \S \ref{lift}) and settled the missing types.





Denote by $\X$ the space of all spatial quartics; it is subdivided into equisingular strata  $\X(S)$ where $S$ is a set of simple singularities. Each stratum $\X(S)$ is further subdivided into its connected components corresponding to equisingular deformation families. The equisingular strata $\mathcal{X}(S)$ splits also into families $\mathcal{X}_*(S)$ where the subscript $*$ is the sequence of invariant factors of a certain finite group $\K$ (see \S \ref{conf} )
A complete description of the strata $\mathcal{X}_1(S)$ constituted by the so-called \emph{nonspecial} quartics, \emph{i.e.}, $\K=0$, is given by G\"{u}ne\c{s} Akta\c{s} \cite{Cisem1}. 
We complete this project in this paper by giving a complete description of the whole equisingular strata $\mathcal{X}(S)$ consisting also the quartics with $\K \ne 0$.
 
A standard \emph{real structure} (\emph{i.e.}, an antiholomorphic involution) $\operatorname{conj}\colon \mathbb{P}^3\rightarrow\mathbb{P}^3$ induces a real structure $\operatorname{c}: \X\rightarrow\X$ which takes a quartic to its conjugate. This real structure $\operatorname{c}$ depends on the choice of $\operatorname{conj}$; however, 
the induced action on the set of the connected components of the equisingular strata is well defined. A connected component $\mathcal{D}\subset\X(S)$ is called \textit{real} if $\operatorname{c}(\mathcal{D})=\mathcal{D}$. 
Clearly, each stratum $\X(S)$ consists of real and pairs of complex conjugate components; this classification of components is given in~\cite{Alex2} for sextics  and in~\cite{Cisem1} for nonspecial quartics.

By a \emph{perturbation} of a set of singularities $S$, we mean a primitive root sublattice $S'$ of the corresponding root lattice $S$ (see \S\ref*{perturbations}). Recall that  unlike high degress,  for a simple quartic surface $X$ with the set of singularities $S$, any perturbation of $S$ is actually realized by a perturbation of $X$. Our extremal families are extremal in the sense that they are not obtained by perturbation  from any bigger family. The list of all equisingular strata $\X(S)$ is way too huge to be 
itemized explicitly; thus in the existence part of the statement of Theorem \ref{th.classification.quartics}, we describe only those strata that are extremal in terms of formal perturbations (with the imprimitivity taken into account). 
 \begin{table}
 \caption{Quartic surfaces with simple singularities}\label{table:quartics}
 \vbox{\def\sep{\ \,}%
 \def\b#1{\setbox0\hbox{$00$}\hbox to\wd0{\hss$#1$\hss}}%
 \let\0\relax
 \offinterlineskip
 \halign{\strut\sep\hss#\hss\sep\sep&&\hss$#$\hss\sep\cr
 \noalign{\hrule}
 $\mu$&\b1&\b2&\b3&\b4&\b5&\b6&\b7&\b8&\b9&10&11&12&13&14&15&16&17&18&19&\text{Total}\cr
 \noalign{\hrule}
 ss&1&2&3&6&9&16&24&39&57&88&128&193&276&403&563&765&880&738&278&4469\cr
 cf&1&2&3&6&9&17&26&46&74&130&211&361&580&939&1370&1779&1766&1178&347&8845\cr
 $r$&1&2&3&6&9&17&26&46&74&130&211&361&580&939&1370&1778&1765&1167&304&8789\cr
 $c$& & & & & & & & & & & & & & & &1&1&11&86&99\cr
 ex& & & & & & & & & & & & & & & &2&1&36&\text{all}&39\cr
 \noalign{\hrule}
 }}
 \end{table}

 The ultimate result can be stated as follows.
 
 \theorem\label{th.classification.quartics}
 The equisingular deformation families of simple quartic surfaces $X\subset\Cp3$ are
 the perturbations of the $390$ maximizing \rom($\mu=19$\rom) families listed in \autoref{table:maxtable}  and $39$
 \emph{extremal} families listed in \autoref{table:extremal39}.
 
 In particular, a non-maximizing quartic is uniquely determined by its configuration (see Definition \ref{configuration.admisibility}) up to deformation and complex conjugation, and there are only $13$ non-maximizing families listed in \autoref{table:nonreal} that are not real.
 \endtheorem
 
The counts are summarized in \autoref{table:quartics} where we list the following data,
itemized by the total Milnor number $\mu:=\rank S$:
\roster*
\item
ss $=$ the number of sets of singularities;
\item
cf $=$ the number of \emph{configurations}(see Definition \ref{configuration.admisibility}), or
\emph{lattice types} in \cite{Shimada:Zsplitting};
\item
$(r,c)$ $=$ the numbers of real components and pairs of complex conjugate
ones;
\item
ex $=$ the number of extremal families.
\endroster
 




 
\subsection{Contents of the paper}
In \S\ref{integral.lattices}, after recalling basic notions and facts of Nikulin's \cite{Niku2} theory of discriminant forms and lattice extensions, we give a brief introduction to Miranda--Morrison's theory \cite{MM1, MM2,MM3} and recast some of their results in a more convenient way to use for our purposes. Then in \S\ref{reflections.and.lifts}, as one of the principal novelties of this paper, we introduce an approach for lifting reflections from finite quadratic form to the $p$-adic lattices which resolves exceptional remaining cases in our computations.

In \S\ref{section.quartics}, we consider quartics as $K3$-surfaces, define the notions of configurations and their realizations and use the theory of $K3$-surfaces to reduce the original geometric problem to a purely arithmetical question concerning realizations of configurations. 

Having the conceptual part settled, the principal result of the paper stated in \S\ref{principal.results} is proved in \S\ref{algorithm.for.classification} where we outline the algorithm used to enumerate the equisingular deformation classes of quartics and give some examples illustrating the steps listed in the general scheme. Finally in \S\ref{example.for.proof}, we demonstrate the proof of Theorem\ref{th.classification.quartics} in a detailed way for the particular set of singularities $S=\mathbf{A}_{15}\oplus\mathbf{A}_3$ and then make some concluding remarks. The tables referred in the main result Theorem \ref{th.classification.quartics} are given in \S\ref{tables}.

\subsection{Acknowledgements}
I am grateful to Alexander Degtyarev for a number of comments, suggestions and fruitful and motivating discussions.


\section{Integral lattices}\label{integral.lattices}

In this introductory section we recall briefly a few elementary facts concerning integral lattices, their discriminant forms and extensions. The principal reference is \cite{Niku2}.

\subsection{Finite quadratic forms (see \cite{MM1,Niku2})}

A \emph{finite quadratic form} is a finite abelian group $\mathcal{L}$ equipped with a map $q\colon \mathcal{L}\rightarrow\mathbb{Q}/2\mathbb{Z}$ quadratic in the sense that 
 $$q(x+y)=q(x)+q(y)-2b(x,y),\quad q(nx)=n^2q(x),\quad x,y\in\mathcal{L},\;n\in \Z,$$ 
 where $b\colon \mathcal{L}\otimes\mathcal{L}\rightarrow\mathbb{Q}/\mathbb{Z}$ is a symmetric bilinear form (which is determined by $q$) and $2\colon \mathbb{Q}/\mathbb{Z}\rightarrow\mathbb{Q}/2\mathbb{Z}$ is the natural isomorphism. We abriviate  $x^2:=q(x)$ and $x\cdot y:=b(x,y)$. 
 
Each finite quadratic form can be decomposed into the orthogonal direct sum $\mathcal{L}=\bigoplus_p\mathcal{L}_{p}$ of its $p$-primary components $\mathcal{L}_{p}:=\mathcal{L}\otimes \mathbb{Z}_p$, where the summation runs over all primes $p$. The \emph{length} $\ell( \mathcal{L})$ is the minimal  number of generators of $\mathcal{L}$; we put $\ell_p(\mathcal{L}):=\ell(\mathcal{L}_{p})$.
A finite quadratic form $\mathcal{L}$ is called \emph{even} if $x^2=0\bmod \Z$ for each element $x\in \mathcal{L}_{2}$ of order $2$; it is called \emph{odd} otherwise.

A finite quadratic form is \emph{nondegenerate} if the homomorphism
\begin{align*}
\mathcal{L}\rightarrow \operatorname{Hom}(\mathcal{L},\mathbb{Q}/\mathbb{Z}),\quad x\mapsto(y\mapsto x\cdot y)
\end{align*}
is an isomorphism. We denote by $\Aut (\mathcal{L})$ the group of automorphisms of $\mathcal{L}$ preserving the form $q$. A subgroup $\K\subset \mathcal{L}$ is called \emph{isotropic} if the restriction of the quadratic form $q$ on $\mathcal{L}$  to $\K$ is identically zero. If this is case $\K^{\bot}/\K$ also inherits from $\mathcal{L}$ a nondegenerate quadratic form.

For a fraction $\frac{m}{n}\in \Q/ 2\Z $, with $(m,n)=1$  such that $mn=0\bmod 2$, we denote by $\langle\frac{m}{n}\rangle$  the nondegenerate finite quadratic form on $\Z/n\Z$ sending the generator to $\frac{m}{n}$, \ie , $\alpha^2=\frac{m}{n} \bmod 2\Z$ for a generator $\alpha$. 
For an integer $k\ge 1$, let $\mathcal{U}_n$ and $\mathcal{V}_n$ be the length $2$ forms on $(\Z/n\Z)^2$, defined by the matrices
\begin{align*}
\mathcal{U}_n:=\left[
\begin{array}{cc}
0 & 1/n \\
1/n & 0 
\end{array}
\right],\quad \mathcal{ V}_n:=\left[
\begin{array}{cc}
2/n&  1/n\\
1/n & 2/n
\end{array}\right], \quad \mbox{where $n=2^k$}.
\end{align*}
Nikulin \cite{Niku2} proved that, a nondegenerate finite quadratic form splits into an orthogonal direct sum of cyclic forms $\langle \frac{m}{n}\rangle$ (defined on the cyclic group  $\Z/n\Z$) and length 2 blocks $\mathcal{U}_n$, $\mathcal{V}_n$. Unless the $2$-torsion  consists of the summands of length 2, we describe nondegenerate finite quadratic forms by expressions of the form $\langle q_1\rangle\ldots \langle q_r\rangle$, where $q_i=\frac{n_i}{m_i}\in \Q$ as above; the group is generated by pairwise orthogonal elements $\alpha_1,\ldots\alpha_n$ (numbered in the order of appearance) so that $\alpha_i^2=\frac{m_i}{n_i} \bmod 2\Z$ and order of $\alpha_i$ is $n_i$.

\begin{definition}\label{detp}
	Let $\mathcal{L}$ be a nondegenerate quadratic form. Given a prime $p$, the determinant of the Gram matrix (in any minimal basis) of $\mathcal{L}_{p}$ has the form  $u/|{\mathcal{L}_{p}}|$ for some unit $u\in\Z_p^{\times}$, and this unit is well defined  modulo $(\Z_p^{\times})^2$ unless $p=2$ and $\mathcal{L}_2$ is odd; in the latter case, $u$ is well defined  modulo the subgroup generated by $(\Z_2^{\times})^2$ and $5$. We define $\det_p \mathcal{L}=u/|{\mathcal{L}_p}|$ where $u\in\Z_p^{\times}/(\Z_p^{\times})^2$ or ${u\in\Z_2^{\times}/(\Z_2^{\times})^2\times\{1,5\}}$ is as above (see \cite{MM3}).
\end{definition}
\begin{remark}\label{NikuDef}
	According to Nikulin~\cite{Niku2}, given a prime $p$ and a quadratic form $\mathcal{L}$ on a $p$ group, there is a $p$-adic lattice $L$ such that $\operatorname{rk} L= \ell_p(\mathcal{L})$ and $\disc L = \mathcal{L}_p$. Unless $p=2$ and $\mathcal{L}$ is odd, such a lattice $L$ is determined by $\mathcal{L}$ uniquely up to isomorphism. In the exceptional case $p=2$ and $\mathcal{L}$ is odd, there are two such lattices that differ by determinants.  One has $\det L= \det_p\mathcal{L}|{\mathcal{L}_p}|^2= u |{\mathcal{L}_p}|$ for  some unit $u$ as in the Definition \ref{detp} (Nikulin uses this equality as a definition of $\det_p \mathcal{L}$).
\end{remark}


\subsection{Integral lattices and discriminant forms}
An \emph{(integral) lattice} is a finitely generated free abelian group $L$ equipped with a symmetric bilinear form $b\colon L\otimes L\rightarrow \mathbb{Z}$. Whenever the form is fixed, we use the abbreviation $x^2:=b(x,x)$ and $x\cdot y:=b(x,y)$. A lattice $L$ is called \emph{even} if $x^2 =0\mod 2$ for all $x\in L$; it is called \emph{odd} otherwise. The \emph{determinant} $\det L \in \Z$ is the determinant of the Gram matrix of $b$ in any integral basis of $L$. 
 A lattice $L$ is called \emph{unimodular} if $\det L=\pm 1$; it is called \emph{nondegenerate} if $\det L \neq 0$, or equivalently, the \emph{kernel}
\begin{align*}
\operatorname{ker}L=L^{\bot} :=\{x\in L \mid\text{ $x\cdot y= 0$ for all $y\in L$}\}
\end{align*}
is trivial.

Given a lattice $L$, the bilinear form on $L$ can be extended by linearity to a $\mathbb{Q}$-valued bilinear form on $L\otimes\mathbb{Q}$. The inertia indexes $\sigma_{\pm}$ and the signature $\sigma:=(\sigma_-,\sigma_+)$ of $L$ are defined as those of $L\otimes\mathbb{Q}$. A nondegenerate lattice is called  \emph{hyperbolic} if $\sigma_+L=1$. 
If $L$ is nondegenerate, then we have  canonical inclusion
\begin{align*}
L\subset L^{\vee}:=\operatorname {Hom}(L,\Z)=\{x \in L\otimes\Q \mid \text{$x\cdot y \in \Z$ for all $y \in L$}\}
\end{align*}
The finite quotient group $\discr L :=L^{\vee}/L$ of order $\mathopen|{\det L}\mathclose|$ is called the \emph{discriminant group}  of $L$. In particular, $L$ is unimodular if and only if $\discr L=0$, \emph{i.e.}, $L=L^{\vee}$

The discriminant group inherits from $L\otimes\mathbb{Q}$ a nondegenerate symmetric bilinear form
\begin{align*}
b\colon \discr L\otimes \discr L\rightarrow \mathbb{Q}/\mathbb{Z},\quad (x\bmod L)\otimes(y\bmod L)\mapsto(x\cdot y)\bmod \Z,
\end{align*}
and if $L$ is even, its quadratic extension
\begin{align*}
q\colon \discr L\rightarrow \mathbb{Q}/2\Z, \quad (x \bmod L) \mapsto x^2 \bmod2\Z,
\end{align*}
called, respectively,  the \emph{discriminant bilinear form} and \emph{discriminant quadratic form}. Note that the discriminant group of an even lattice is a finite quadratic form. We use the notation $\discr_p L$ for the $p$-primary part of $\discr L$.
When speaking about the discriminant groups and  their (anti-)isomorphisms, these forms are always taken into account.

Lattices are naturally grouped into \emph{genera}. Omitting the precise definition we follow Nikulin \cite{Niku2} who states that two nondegenerate even lattices $L',L''$ are in the same \emph{genus} if and only if one has $\operatorname{rk} L'=\operatorname{rk} L''$, $\sigma L'=\sigma L''$ and $\discr L'\cong\discr L''$. Each genus consists of finitely many isomorphism classes.

An isometry  $\psi \colon L \rightarrow L'$ between two lattices is a group homomorphism respecting the bilinear forms; obviously one always has $\operatorname{Ker} \psi \subset \operatorname{Ker} L$. The group of bijective autoisometries of a nondegenerate lattice $L$ is denoted by $O(L)$. The action of $O(L)$ extends to $L\otimes\mathbb{Q}$ by linearity, and the latter action descents to  $\discr L$. Therefore, there is a natural homomorphism $O(L)\rightarrow \Aut(\discr L)$ where $\Aut(\disc L)$ denotes the group of automorphisms of $\discr L$ preserving the discriminant form $q$ on $\discr L$. In general this map is neither one-to-one nor onto; however, 
without any confusion we freely apply autoisometries $g\in O(L)$ to objects in $\discr L$. Obviously one has $\Aut(\discr L)=\prod_p \Aut(\discr_p L)$ where the product runs over all primes. The restriction of $d$ to $p$-primary components are denoted by $d_p\colon O(L)\rightarrow \Aut(\discr_p L)$.

A $4$-\emph{polarized lattice} is a nondegenerate hyperbolic lattice $L$ equipped with a distinguished vector $h\in L$ such that $h^2=4$; this vector is usually assumed but often omitted from the notation. The group of polarized autoisometries is denoted by $O_h(L)$, \emph{i.e.}, the group of automorphisms of $L$ preserving $h$. 

The orthogonal projection establishes a liner isomorphism between any two maximal positive definite  subspaces in $L\otimes\mathbb{R}$, thus providing a way for comparing orientations. A coherent choice of orientations of all maximal positive definite subspaces is called a \textit{positive sign structure} on $L$. We denote by $O^+(L)\subset O(L)$ the subgroup consisting of the autoisometries preserving a positive sign structure. Either one has $O^{+}(L)=O(L)$ or $O(L)^+$ is a subgroup of $O(L)$ of index $2$. In the latter case, each element of $O(L)\smallsetminus O^+(L)$ is called a \emph{skew-autoisometry} of L,\emph{ i.e.}, the autoisometries of $L$ that reverse the positive sign structure.

Of special importance are so called reflections of $L$: given a nonzero vector $a\in L$, the \emph{reflection} defined by $a$ is the automorphism 
\[t_a\colon L\rightarrow L,\quad x\mapsto \frac{2(a\cdot x)}{a^2}a .\]
It is well defined if and only if $(2a/a^2)\in L^\vee$. Note that $t_a$ is an involution. Each image $d_p(t_a)$ is also reflection and if $a^2=\pm1$ or $a^2=\pm2$, then the induced automorphism $d(t_a)$ of the discriminant group is the identity and $t_a$ extends to any lattice containing $L$.

\subsection{Root Lattices}\label{root.lattices}

A \emph{root} in an even lattice $L$ is a vector $r\in L$ of square $(-2)$. A \emph{root lattice} is an even negative definite lattice generated by its roots. Each root lattice has a unique decomposition into orthogonal direct sum of irreducible root lattices, the latter being those of types $\textbf{A}_n$, $n\geq1$, $\textbf{D}_n$, $n\geq4$, or $\textbf{E}_n$, $n=6,7,8$. 

Given a root lattice $S$, the vertices of the Dynkin diagram $\Gamma:=\Gamma_S$ can be identified with the elements of a basis for $S$ constituting a single Weyl chamber. Thus, one has an obvious homomorphism $\operatorname{Sym}(\Gamma)\rightarrow O(S) $
 where $\operatorname{Sym}(\Gamma)$ is the group of  symmetries of the Dynkin diagram $\Gamma$. By the classification of the connected Dynkin graphs, for irreducible root lattices, the groups $\operatorname{Sym}(\Gamma)$ are given as follows:
 \begin{enumerate}
 	\item if $S$ is $\mathbf{A}_1$, $\mathbf{E}_7$ or $\mathbf{E}_8$, then $\operatorname{Sym}(\Gamma)=1$
 	\item if $S$ is $\mathbf{D}_4$, then $\operatorname{Sym}(\Gamma)=\mathbb{S}_3$
 	\item for all other types, $\operatorname{Sym}(\Gamma)=\mathbb{Z}_2$
 \end{enumerate}

If S is $\mathbf{A}_p$, $p\geq2$, $\mathbf{D}_{2k+1}$ or $\mathbf{E}_8$, then the only nontrivial symmetry of $\Gamma$ induces $-\operatorname{id}$ on $\discr S$. If $S$ is $\mathbf{E}_8$ then $\discr S=0$ and if $S$ is $\mathbf{A}_1$, $\mathbf{A}_7$ of $\mathbf{D}_{2k}$, the groups $\discr S$ are $\mathbb{F}_2$ modules  and $-\operatorname{id}=\operatorname{id}$ on $\Aut (\discr S)$.


\subsection{Lattice extensions}\label{lattice.extensions}

From now on unless specified otherwise all lattices considered are even and nondegenerate. An extension of an even lattice $S$ is an even lattice $L$ containing $S$. Two extensions $L',L''\supset S$ are called \emph{isomorphic} if there is a bijective isometry $L'\rightarrow L''$ preserving $S$, in particular, if the isomorphism  $L'\rightarrow L''$ is identical on $S$, the extensions $L'$ and $L''$ are called \emph{strictly isomorphic}. More generally, one can also fix a subgroup $G\in O(S)$ and speak about $G$-\emph{isomorphisms} of the extension, \emph{i.e.}, bijective isometries whose restriction to $S$ is in $G$.

The two extreme cases are \emph{finite index extensions}, \emph{i.e.}, $L$ contains $S$ as a subgroup of finite index and \emph{primitive extensions}, \emph{i.e.}, $L/S$ is torsion free. The general case $L\supset S$ splits into the finite index extension $\tilde{S}\supset S$ and primitive extension $L\supset\tilde{S}$, where 
$$\tilde{S}:=\{x\in L\mid nx\in S\;\mbox{for some $n\in\Z$}\}$$
is the \emph{primitive hull} of $S$ in $L$. 

Any extension $L\supset S$ of finite index admits a unique embedding $L\subset S\otimes \Q$. Since $S$ is nondegenerate, we have $L\subset S^\vee$, and hence the natural inclusions
\begin{align*}
S\subset L\subset L^{\vee}\subset S^{\vee}.
\end{align*}  
The subgroup $\K:=L/S\subset S^{\vee}/S=\discr S$ is called the \emph{kernel} of the finite index extension $L\supset S$. This subgroup $\K$ is \emph{isotropic} (since $L$ is an even integral lattice), \emph{i.e.}, the restriction to $\K$ of the quadratic form $q\colon \discr S\rightarrow \mathbb{Q}/2\Z $ is identically zero. Conversely, if $\K\subset \discr S$ is isotropic, the lattice
 $$L:=\{x\in S\otimes \mathbb{Q} \mid x\bmod S\in \K\}$$
  is an extension of $S$ and we say that $L$ is the extension of $S$ by $\K$. Thus, we have the follwing statement. 
\begin{theorem}[Nikulin \cite{Niku2}]\label{L-K}
	Let $S$ be a nondegenerate even lattice, and fix a subgroup $G\subset
	O(S)$. The map
	\begin{align*}
L\mapsto \mathcal{K}:=L/S \subset \discr S
	\end{align*} 
	establishes a one-to-one correspondence between the set of
	$G$-isomorphism classes of finite index extensions $L\supset S$
	and the set of $G$-orbits of isotropic subgroups
	$\mathcal{K}\subset \discr S$. Under this correspondence one
	has $\discr L=\mathcal{K}^{\bot}/\mathcal{K}$. Furthermore
	\begin{enumerate}
		\item an autoisometry $g\in O(S)$ extends to  $L$ if and only if $g(\K)=\K$;
		\item two extensions $L'$, $L''$ are isomorphic if and only if their kernels $\K'$, $\K''$ are in the same $O(S)$-orbit, \emph{i.e}, there is $g\in O(S)$ such that $g(\K')=g(\K'')$.
	\end{enumerate}
	 
\end{theorem}
Another extreme case is that of a primitive extension $L\supset S$, (\emph{i.e.}, such that the group $L/S$ is torsion free). Such extensions are studied by fixing (the isomorphism class of) the orthogonal complement $T:=S^{\bot}\in L$. Then $L$ is a finite index extension of $S\oplus T$, in which $T$ is also primitive, and by Theorem \ref{L-K}, it is described by an isotropic subgroup
\begin{equation*}
\K\subset \discr (S\oplus T)=\discr S\oplus \discr T.
\end{equation*}
Since $S$ and $T$ are both primitive in $L$, the kernel $\mathcal{K}$  does not intersect with any of $\discr S$ and $\discr T$. It follows that the projection maps
\begin{equation*}
\operatorname{proj}_{S}\colon \K \rightarrow \discr S\mbox{ and } \operatorname{proj}_{T}\colon \K \rightarrow \discr T
\end{equation*}
are both monomorphisms. Since $\K$ is isotropic, it is the graph of a bijective anti-isometry $\psi\colon \discr S'\rightarrow\discr T'$, where $\discr S'= \operatorname{proj}_{S}(\K)$ and $\discr T'= \operatorname{proj}_{T}(\K)$. Conversely, given a bijective anti-isometry $\psi\colon\discr S'\rightarrow\discr T'$ where $\discr S'\subset\discr S$ and $\discr T'\subset \discr T$, the graph of $\psi$ is an isotropic subgroup $\mathcal{K}\subset\discr S\oplus \discr T$ and the corresponding finite index extension $L\supset S\oplus T$ is a primitive extension whose kernel is $\K$. Thus, we have the following statement (\cf. Nikulin~\cite{Niku2}).
\begin{lemma}\label{genlemma}
	Given two nondegenerate even lattices $S$, $T$ and a subgroup $G\subset
	O(S)\times O(T)$, there is a one-to-one correspondence between the set of $G$-isomorphism classes of finite index extensions $L\supset S\oplus T$  in which both $S$ and $T$ are primitive  and that of $G$-conjugacy classes of bijective anti-isometries
	\begin{equation}\label{isodisc}
	\psi\colon \discr S'\rightarrow \discr T'
	\end{equation}
	where $\discr S'\subset \discr S$ and $\discr T'\subset \discr T$. Furthermore, a pair of isometries $f\in O(S)$ and $g\in O(T)$ extends to $L$ if and only if $f|_{\discr S'}=\psi^{-1}g|_{\discr T'}\psi$ in $\Aut (\discr S')$.
\end{lemma}
If $L$ above is unimodular, $\disc L=0$, we have
$\mathopen|\discr S\mathclose|\mathopen|\discr T\mathclose|=\mathopen|\discr S'\mathclose|\mathopen|\discr T'\mathclose|$. Hence, $\discr S'=\discr S$ and $\discr T'=\discr T$ and $\psi$ in \eqref{isodisc} is an anti-isomorphism $\discr S\rightarrow\discr T$. Since also $\sigma_{\pm}T=\sigma_{\pm}L-\sigma_{\pm}S$, it follows that the genus $g(T)$ is determined by the genera $g(S)$ and $g(L)$; we will denote this common genus by $g(S^{\bot}_L)$ (We emphasize that $g(S^{\bot}_L)$  merely encodes a ``local data" composed formally from $g(S)$ and $g(L)$; apriori, it may even be empty, \emph{cf.} Theorem \ref{th.N.existence} below). If $L$ is also indefinite, it is unique in its genus (see, \emph{e.g.}, Siegel~\cite{SiegelI,SiegelII,SiegelIII}). Then we have the following corollary of the above lemma.
\begin{corollary}\label{bicoset}
	Given  a subgroup $G\subset O(S)$ and unimodular even indefinite lattice $L$, a primitive isometry  $S \into L$ give rise to a bijective isometry $\psi\colon \discr S\rightarrow -\discr S^\bot$ and $G$-isomorphism classes of a primitive isometries  $S\into L$ are in canonical bijection with the following sets of data:
	\begin{enumerate}
		\item an even lattice (isomorphism class) $T\in g(S^{\bot}_L)$, and
		\item a bi-coset in $G\backslash \Aut (\discr T)/O(T)$.
	\end{enumerate}
\end{corollary}

In particular, the extension $L\supset S$ exists if and only if the genus $g(S^{\bot}_L)$ is nonempty.

From now on, we fix the notation $\L:=3\textbf{U}\oplus2\textbf{E}_8$ where $\textbf{U}$ stands for the \emph{hyperbolic plane}, the lattice generated by a pair of vectors $u,v$ (refered as the \emph{standard basis} of $\textbf{U}$) with $u^2=v^2=0$ and $u\cdot v=1$. Note that $3\textbf{U}\oplus2\textbf{E}_8$ is the unique even unimodular lattice of rank $22$ and signature $(3,19)$. We are concerned about this lattice as it is the intersection index form of a $K3$-surface $X$, \ie, $H_2(X;\Z)\cong \L$.   
We are interested in the primitive embeddings to $\L$. The following theorem giving a criterion for $g(S^{\bot}_\textbf{L})\neq\emptyset$ is a combination of the above observation and Nikulin's existence theorem~\cite{Niku2} applied to the genus $g(T)$.
\begin{theorem}[Nikulin~\cite{Niku2}]\label{th.N.existence}
	Given a nondegenerate even lattice~$S$, a primitive extension $\textbf{L}\supset S$ exists
	if and only if
	the following conditions hold
	\begin{enumerate}
		\item $\sigma_+ S\leq 3$, $\sigma_-S\leq 19$ and $\ell(\mathcal{S})\leq 22- \operatorname{rk} S$, where $\mathcal{S}=\discr S$;
		\item one has  $|{\mathcal{S}}|\det_p (\mathcal{S}) =(-1)^{\sigma_+S-1} \bmod (\Z_p^{\times})^2$ for each  odd prime $p$ such that  $\ell_p(\mathcal{S})= 22- \operatorname{rk} S$;
		\item If $\ell_2(\mathcal{S})= 22- \operatorname{rk} S$, and $\mathcal{S}_2$ is even  then $|{\mathcal{S}}|\det_2 (\mathcal{S}) =\pm 1 \bmod (\Z_2^{\times})^2$.
	\end{enumerate}
\end{theorem}

\subsection{Miranda-Morrison Theory}\label{MM}
Following the classical approach,  Nikulin \cite{Niku2} gives the sufficient conditions to obtain the uniqueness of an even indefinite lattice $T$ of rank at least $3$ in its genus and surjectivity of the map $d\colon O(T)\rightarrow \Aut (\discr T)$. However, those conditions do not capture all the cases that we want to cover, hence we apply the stronger (non)uniqueness criteria due to Miranda- Morrison \cite{MM1, MM2, MM3} extending Nikulin's work. Throughout this section, we assume that $T$ is an indefinite nondegenerate even lattice of rank $\operatorname{rk} T \geq 3$. 

With the ultimate goal of calculating the groups $E(T)$ and $E^+(T)$ see \eqref{E(N)} and \eqref{E(N)plus}, respectively, and the images of some certain maps in these groups, it is convenient to introduce following groups 
\begin{align*}
\Gamma_p:&=\{\pm1\}\times\mathbb{Q}^{\times}_p/(\mathbb{Q}^{\times}_p)^2,\\
\Gamma_0:&=\{\pm1\}\times\{\pm1\}\subset\Gamma_{\mathbb{Q}}:=\{\pm1\}\times\mathbb{Q}^{\times}/(\mathbb{Q}^{\times})^2,
\end{align*}
and following subgroups related to $\Gamma_p$ :
\begin{itemize}
\item $\Gamma_{p,0}:=\{(1,1),(1,u_p),(-1,1),(-1,u_p)\}\subset\Gamma_p$; here, $p$ is odd and $u_p$ is the only nontrivial element of $\mathbb{Z}^{\times}_p/(\mathbb{Z}^{\times}_p)^2$,
\item $\Gamma_{2,0}:=\{(1,1),(1,3),(1,5),(1,7),(-1,1),(-1,3),(-1,5),(-1,7)\}\subset\Gamma_2$,
\item $\Gamma_{2,2}:=\{(1,1),(1,5)\}\subset\Gamma_2^{++}$,
\item $\Gamma_0^{--}:=\{(1,1),(-1,-1)\}\subset\Gamma_0$.
\end{itemize}
We also define
\begin{align*}
\Gamma_{\mathbb{A},0}:=\prod_p\Gamma_{p,0}\subset\Gamma_{\mathbb{A}}:=\Gamma_{\mathbb{A},0}\cdot\sum_p\Gamma_p\subset \Gamma :=\prod_p\Gamma_p
\end{align*}
where we use $\cdot$ to denote the sum of subgroups, while we reserve the notation $\sum$ and $\prod$ to distinguish between direct sums and products.
Note that
\begin{align*}
\Gamma_{\mathbb{A}}=\{(d_p,s_p)\in \Gamma\ |(d_p,s_p)\in\Gamma_{p,0}\text{ for almost all $p$}\}
\end{align*}
Defined in the notation of \cite{MM3}, we will use certain subgroups
$$\Sigma_p^{\sharp}(T):= \Sigma^{\sharp}(T\otimes\mathbb{Z}_p)\quad\mbox{and}\quad \Sigma_p(T):=\Sigma(T\otimes\mathbb{Z}_p)$$
which are both a priori subgroups of $\Gamma_p$ (we refer the reader to \cite{MM3}, see chapter 7, section 4, for the precise definitions). In fact $\Sigma_p^{\sharp}\subset \Gamma_{p,0}$ always and $\Sigma_p\subset \Gamma_{p,0}$ for almost all p. The subgroups $\Sigma^{\sharp}_p(T)$ are computed explicitly in \cite{MM3} (see Theorems 12.1, 12.2, 12.3 and 12.4 in  chapter 7). One has
$$\Sigma^{\sharp}(T):= \prod_p\Sigma^{\sharp}_p(T)\subset\Gamma_{\mathbb{A},0}\quad\mbox{and}\quad \Sigma(T):= \prod_p\Sigma_p(T)\subset\Gamma_{\mathbb{A}}.$$
We introduce the \emph{Miranda-Morrison group} $E(T)$ as it is defined  in \cite{MM3} (see chapter $8$, sections $5$, $6$ and $7$):
\begin{align}\label{E(N)}
E(T):=\Gamma_{\mathbb{A},0}/\prod_p\Sigma_p^{\sharp}(T)\cdot\Gamma_0.
\end{align}
Crucial is the fact that $\Sigma^{\sharp}_p(T)=\Gamma_{p,0}$ unless $p\mathrel|\operatorname{det}(T)$; thus, (\ref{E(N)}) reduces to finitely many primes $p$:
\begin{align}\label{E(N)finite}
E(T)=\prod_{p|\operatorname{det}(T)}\Gamma_{p,0}/\prod_{p|\operatorname{det}(T)}\Sigma_p^{\sharp}(T)\cdot\Gamma_0.
\end{align}
Hence, this group is finite.  We call a prime $p$ \emph{irregular} with respect to $T$ if $p\mid\operatorname{det}(T)$.

Consider the natural map $\mathbb{Q}^{\times}/(\mathbb{Q}^{\times})^2\rightarrow\mathbb{Q}^{\times}_p/(\mathbb{Q}^{\times}_p)^2$ inducing the projections
\begin{align*}
\varphi_p:\Gamma_0\rightarrow\Gamma_{p,0},
\end{align*}
we define the invariants 
\begin{align*}
e_p(T):=[\Gamma_{p,0}:\Sigma^{\sharp}_p(T)] \text{ and } \tilde{\Sigma}_p(T)=\Sigma^{\sharp}_0(T\otimes\mathbb{Z}_p):=\varphi^{-1}_p(\Sigma^{\sharp}_p(T))\subset\Gamma_0,
\end{align*}
used in the following theorem.
\begin{theorem}[Miranda--Morrison \cite{MM3}]\label{MMexact.sequence}
	Let $T$ be a non-degenerate indefinite even lattice with $\operatorname{rk}(T)\geq 3$. Then there is an exact sequence
	\begin{align}\label{exactsequence}
	O(T)\xrightarrow{d}\operatorname{Aut}(\operatorname{discr} T)\xrightarrow{\mathrm{e}} E(T)\rightarrow g(T)\rightarrow 1,
	\end{align}
	where $g(T)$ is the genus group of $T$. One has
	\begin{align}\label{orderofEN}
	|E(T)|=\frac{e(T)}{[\Gamma_0:\tilde\Sigma(T)]}
	\end{align}
	where
	\begin{align*}
	e(T):=\prod_{p|\operatorname{det}(T)} e_p(N), \quad\tilde\Sigma(T):=\bigcap_{p|\operatorname{det}(T)} \tilde\Sigma_p(T),
	\end{align*}
\end{theorem}
Algorithims computing $e_p(N)$ and $\tilde{\Sigma}_p(N)$ explicitly are given in \cite{MM2}. Computations are in terms of $\operatorname{rk} T$, $\operatorname{det} T$ and $\discr T$ only, it follows that the genus group $g(T)$ determines $e_p(T)$ and $\tilde{\Sigma}_p(T)$ and also $\operatorname{Coker}(d\colon O(T)\rightarrow \Aut (\discr T))$. Computing the Miranda-Morison group $E(T)$ is even easier than computing its constituents,  $\operatorname{Coker}(d)$ and the genus group $g(T)$. Since one can read the the subgroups  $\tilde{\Sigma}_p(T)$ and $\Sigma^{\sharp}_p(T)$ from the tables given in \cite{MM1} (see chapter 7 section 12), the computation is immediate.

As an unimodular even indefinite lattice is unique in its genus, one can obtain the next statement by combining the Corollary \ref{bicoset} and Theorem \ref{MMexact.sequence}.
\begin{theorem}[Miranda--Morrison \cite{MM1,MM2}]\label{correspondence.E(N)}
	Let $S$ be a primitive sublattice of an even unimodular lattice $L$ such that $T:=S^{\bot}$ is a non-degenerate indefinite even lattice with $\operatorname{rk}(T)\geq3$. Then the strict isomorphism classes of primitive extensions $S\hookrightarrow L$ are in a canonical one-to-one correspondence with the group $E(T)$.
\end{theorem}
Given a unimodular lattice $L$ and a primitive sublattice $S\subset L$, fix an anti-isometry $\psi\colon \discr S\rightarrow \discr T$ and consider the induced map $d^{\psi}\colon O(S)\rightarrow\operatorname{Aut}(\discr T)$ (see \S \ref{lattice.extensions}). If $T$ is as in Theorem \ref{correspondence.E(N)}, then $\operatorname{Im} d\subset\Aut(\discr T)$ is a normal subgroup with abelian quotient  and  we have a homomorphism  
\begin{align}\label{dperb}
d^{\bot}\colon O(S)\rightarrow \Aut(\discr T)\xrightarrow{\mathrm{e}}E(T)
\end{align}
independent of $\psi$. Then the next statement generalizing Theorem \ref{correspondence.E(N)} follows from Theorem \ref{MMexact.sequence} and Corollary \ref{bicoset}.
\begin{corollary}\label{coker}
	Let $S$ be a primitive sublattice of an even unimodular lattice $L$ such that $T:=S^{\bot}$ is a non-degenerate indefinite even lattice with $\operatorname{rk}(T)\geq3$ and let $G\subset O(S)$ be a subgroup. Then, the $G$-isomorphism classes of primitive extensions $S\hookrightarrow L$ are in a one-to-one correspondence with the $\mathbb{F}_2$-module $E(T)/d^\bot(G)$.
\end{corollary}
Let $p$ be a prime and consider the homomorphism
\begin{align}\label{product.epimorphisms}
\Aut(\discr T)=\prod_p\Aut(\discr_p T)\xrightarrow{\phi}\prod_p\Sigma_p(T)/\Sigma^{\sharp}_p(T)
\end{align}
which is the product of the epimorphisms
$
\phi_p\colon \Aut(\discr_p T)\twoheadrightarrow\Sigma_p(T)/\Sigma^{\sharp}_p(T).
$ Note that the map
$$\mathrm{e}\colon\prod_p \Aut(\discr T)\rightarrow E(T)$$
given in \eqref{exactsequence} does not preserve product structures. 


\begin{remark}\label{efibeta}
Crucial is the fact that the map $\mathrm{e}$ in \eqref{exactsequence} is given as $\mathrm{e}=\beta\circ\phi$ where $\beta$ is the quotient projection
 $$\beta\colon \prod_p\Sigma_p(T)/\Sigma^{\sharp}_p(T) \rightarrow E(T).$$
\end{remark}

\subsection{Reflections and their lifts}\label{reflections.and.lifts}
Let $p$ be a prime and consider an element $a\in \discr_p T$ satisfying
\begin{equation}\label{star}
\text{$p^ka=0$ and $a^2=\frac{2u}{p^k} \bmod2\mathbb{Z}$, g.c.d$(u,p)=1$, $k\in\mathbb{N}$.}
\end{equation}
Then the map $ x\mapsto2(x\cdot a)/a^2\bmod p^k$ is a well defined functional; thus there is a reflection $t_a\in \Aut(\discr_p T)$,
\begin{align*}
t_{a}\colon x\mapsto x-\frac{2(x\cdot a)}{a^2}a.
\end{align*}
Note that if $2a=0$ 
then $t_a=\operatorname{id}$.

The images of the homomorphism $\phi$ given in \eqref{product.epimorphisms} can easily be computed on the reflections $t_a$ by lifting them to the corresponding $p$-adic lattice: Let $t_{\bar{a}}$ be a lift of $t_a$ to the $p$-adic lattice where $\bar{a}\in T\otimes\Z_p$ such that $a=(\bar a/p^k)\bmod T\otimes\Z_p$.  Then, the image of $\phi_p$ on $t_a$ is given via \emph{spinor norm} $\operatorname{spin}(t_a)$ introduced in \cite{MM3} (see chapter $1.10$) which is essentially defined as
\begin{align}\label{spinn}
 \operatorname{spin}(t_a)=\frac{1}{2}\bar{a}^2 \bmod (\Z_p^{\times})^2.
\end{align}
The map $\phi_p$ is then defined as the spinor norm  modulo the indeterminacy subgroup $\Sigma^{\sharp}_p(T)$, hence, as given in \cite{Alex2}, one has
\begin{align}
\phi_p\colon \Aut(\discr_p T)&\twoheadrightarrow\Sigma_p(T)/\Sigma^{\sharp}_p(T)\nonumber\\
t_a&\mapsto (-1,up^k)\label{imageofref}
\end{align}
where $a$ is as in \eqref{star}. 

To attain the goal of this paper, we need to compute the image $\phi_p(t_a)$, however the main problem is that in our computations we don't know the lattice $T$, we only know its genus by Nikulin \cite{Niku2}, hence our goal is to do this computation in terms of genus, \ie, signature $(\sigma_+ T,\sigma_-T)$ and finite quadratic form  $\mathcal{T}=\discr T$, only. With few exceptions, while computing $\phi_p(t_a)$ as in \eqref{imageofref}, the value of $\bar{a}/2$ is almost always welldefined whenever $u$ is sufficiently welldefined; the exceptions are treated in \S \ref{lift}.

\subsubsection{Principal novelty of the paper}\label{lift}
While computing the image $\phi_p(t_a)$ given in \eqref{imageofref}, the following two cases are exceptional and need special treatment:  
\begin{enumerate}
	\item $p=2$ and $a^2=0\bmod \mathbb{Z}$,
	\item $p=2$ and $a^2=\frac{1}{2}\bmod \mathbb{Z}$,
\end{enumerate}
In these two cases, 
the approach in computing the image $\phi_2(t_a)$ as in \eqref{imageofref} has a disambiguation to be clarified: 
In  \eqref{star}, the value of $u$ in the numerator is defined $\mod\,2\Z$ for the former case and $\mod\,4\Z$ for the latter; which is not enough as $\mod\,(\Z_2^{\times})^2$ is essentially $\mod\,8$ (see \eqref{spinn}). Hence, we need to actually write down the $2$-adic lattice $T\oplus\Z_2$ and actually find a lift $\bar a$. To recover $T\oplus\Z_2$ precisely, we use the partial normal form decomposition given in \cite{MM3} (see Chapter 4.4) for $\discr_2 T$. Basically, to lift each summand in this decomposition one by one,  for each standard rational matrix, we write a standard $2$-adic matrix. 
 
However, there is another problem to be fixed when the discriminant $\discr_2 T$ is odd: If $\discr_2 T$ is odd, according to Nikulin \cite{Niku2}, up to isomorphism there are two different $2$-adic lattices $T'$ and  $T''$ with $\operatorname{rk}(T')=\operatorname{rk}(T'')=\ell_2(\discr T)$ and $\discr T'=\discr T''=\discr_2 T$, the ratio of their determinants being $5\in \Z_2^{\times}/(\Z_2^{\times})^2$, see Remark \ref{NikuDef}. Roughly speaking, if the $\discr_2 T$ contains summands of the form $\langle\pm\frac{1}{2}\rangle$, then each of them lifts  to either $[\pm2]$ or $[\pm10]$. We choose lifts arbitrarily in all such summands  but one and this one remaining lift is adjusted by using $\det T$ which is an invariant of the genus group $g(T)$: Recall that $\operatorname{det} T= (-1)^{\sigma_-}|\discr T|$, thus the correct $2$-adict lift $T^*$ is the one satisfying $\operatorname{det} T^*=(-1)^{\sigma_-}|\discr T|\bmod (\Z_2^{\times})^2$.

Once the correct lift is chosen, to compute the image of $\phi_2(t_a)$,
we expand the vector $a\in \discr_2 T$ in the basis vectors of discriminant group, take a particular lift $\bar a\in T\otimes \Z_2$ with the same coordinate vector and compute $\bar{a}^2$ honestly without reducing it $\mod\,2\Z$  
to finally compute the value 
$\bar{a}^2 \bmod (\Z_2^{\times})^2$
given in \eqref{spinn}


Then, as given in Remark \ref{efibeta}, the images of the map $\mathrm{e}$  on the reflections $t_a\in\Aut(\discr T)$ can be computed via $\phi$ and $\beta$.
 
Defined and computed in \cite{MM2}, we introduce the group
\begin{align}\label{E(N)plus}
E^+(T):=\Gamma_{\mathbb{A},0}/\prod_p\Sigma_p^{\sharp}(T)\cdot\Gamma_0^{--}.
\end{align}
where the product reduces to finitely many primes $p|\operatorname{det} T$ as in \eqref{E(N)finite}.
Then one has an exact sequence
\begin{align*}
O^+(T)\xrightarrow{d}\operatorname{Aut}(\operatorname{discr} T)\xrightarrow{\mathrm{e}^+} E^+(T)\rightarrow g(T)\rightarrow 1,
\end{align*}
as in Theorem \ref{MMexact.sequence}, where for the order $|E^+(T)|$ one only replaces $[\Gamma_0:\tilde\Sigma(T)]$ in (\ref{orderofEN}) with $[\Gamma_0^{--}:\tilde{\Sigma}(T)\cap\Gamma_0^{--}]$.  
Under this setting, as in \eqref{dperb}, we have a well defined homomorphism 
\begin{align}\label{dperbplus}
d^{\bot}_+\colon O(S)\rightarrow \Aut(\discr T)\xrightarrow{\mathrm{e}^+}E^+(T)
\end{align}
Then for an element $a\in \discr_p T$ satisfying \eqref{star}, the image of the map $e^+$  on the reflection $t_a\in\Aut(\discr_p T)$ is given by $\mathrm{e^+}=\beta^+\circ\phi^+$, where the maps $\phi$ and $\beta$ in  Remark \ref{efibeta} are replaced with $\phi^+$ and $\beta^+$ by replacing the group  $\Gamma_0$ and any subgroup $H$ of it with $\Gamma_0^{--}$ and $H\cap \Gamma_0^{--} $, respectively. 


\begin{remark}
What is explained in \S\ref{lift}, is not a \emph{panacea} because the group $\Aut(\discr_2 T)$ is not always generated by reflections.
However, experimentally it turns out that  lifting just reflections in $\Aut(\discr_2 T)$ is enough to cover all our needs, see \S \ref{algorithm.for.classification} and furthermore, it appears that it would cover most of $K3$-related problems.
\end{remark}
\begin{remark}
Given in \cite{MM3}, the spinor norm  is computed in terms of reflections, \ie, 
$$\operatorname{spin}(\tau)=\prod \nu_i^2 \bmod (\Z_2^{\times})^2$$
for an element $\tau \in O(T)$ suth that $\tau=t_{\nu_1}t_{\nu_2}\cdots t_{\nu_r}$ where $t_{\nu_i}$ is a reflection against $\nu_i\in T$.\\ 
Shimada \cite{Shimada:Zsplitting} uses an alternative approach: Instead of  decomposing an  automorphism of $\discr_p T$ into reflections and lifting them one by one, he lifts (rather approximates) the whole automorphism and then decomposes it into reflections.
\end{remark}

\section{$K3$-surfaces}\label{section.quartics}
In this section we give a brief introduction to the theory of $K3$-surfaces; for further details, we refer the interested reader to \cite{Huybrechts:K3}. Then we discuss simple quartics as $K3$-surfaces.

\subsection{K3-surfaces}
 A $K3$-surface over $\C$ is a simply connected, compact, complex surface of dimension $2$ whose canonical bundle is trivial. All $K3$-surfaces are K\"{a}hler and since any smooth complete surface is projective, $K3$-surfaces are all projective.

We identify a $K3$-surface $X$ with its polarized N\'{e}ron Severi lattice $NS(X)\ni h$. As is well known, the N\'{e}ron Severi lattice $NS(X)$ of any projective $K3$-surface $X$ is hyperbolic, \emph{i.e.}, $\sigma_+ NS(X)=1$ and admits a primitive embedding 
$$ NS(X)\subset H_2(X;\Z)\cong \mathbf{L}=3\textbf{U}\oplus 2\textbf{E}_8;$$
where $\mathbf{L}$ is the only (up to isomorphism) even unimodular lattice with $\sigma_+\mathbf{L}=3$ and $\sigma_-\mathbf{L}=19$, see \S\ref{lattice.extensions}, hence $\operatorname{rk}NS(X)\le 20$.

\subsection{Quartics as $K3$-surfaces}\label{quartics.as.K3}

A \textit{quartic} is a surface $X\subset\mathbb{P}^3$ of degree four. A quartic is \textit{simple} if all its singular points are simple, i.e., those of type $\mathbf{A},\mathbf{D},\mathbf{E}$, see \cite{Dufree}. 
Given a simple quartic $X\subset\mathbb{P}^3$, its minimal resolution of singularities $\tilde{X}$ is a $K3$-surface; hence, $H_2(\tilde{X})\cong2\mathbf{E}_8\oplus3\mathbf{U}$. We fix the notation $\mathbf{L}_X:=H_2(\tilde{X})$.

For each simple singular point $p$ of $X$, the components of the exceptional divisor over $p$ span a root lattice in $\mathbf{L}_X$. The orthogonal sum of these sublattices, denoted by $S_X$, is identified with the set of singularities of $X$. Recall that the types of individual singular points are uniquely recovered from $S_X$, see \ref{root.lattices}. 

In what follows we identify homology and cohomology of $\tilde{X}$ \emph{via} Poicar\'{e} duality and introduce the following vectors and sublattices:
\begin{itemize}
	\item $S_X\subset \mathbf{L}_X$: the sublattice generated the set of classes of exceptional divisors appearing in the blow-up map $\tilde{X}\rightarrow X$;
	\item $h_X\in \mathbf{L}_X$: the pull-back of the hyperplane section class in $H_2(\mathbb{P}^2)$;
	\item $S_{X,h}=S_X\oplus \mathbb{Z}h_X\subset \mathbf{L}_X$;
	\item $\tilde{S}_X:=(S_X\otimes\mathbb{Q})\cap \mathbf{L}_X$ and $\tilde{S}_{X,h}:=(S_{X,h}\otimes\mathbb{Q})\cap \mathbf{L}_X$ : the primitive hulls of $S_X$ and $S_{X,h}$, respectively; we have $\tilde{S}_X \subset \tilde{S}_{X,h}\subset \mathbf{L}_X$;
	\item $\omega_X \subset \mathbf{L}_X\otimes \mathbb{R}$: the oriented $2$-subspace spanned by the real and imaginary parts of the class of a holomorphic $2$-form on $\tilde{X}$ (the \emph{period} of $\tilde{X}$).
\end{itemize}
Note that $\omega_X$ is positive definite and orthogonal to $h_X$; furthermore, the Picard group $\operatorname{Pic} \tilde{X}$ can be identified with the lattice $\omega_X^\bot\cap \mathbf{L}_X$. In particular $\omega_X\in \tilde{S}_X^\bot\otimes \R$.
The rank $\operatorname{rk}(S_X)$ equals the total Milnor number $\mu(X)$. Since $S_X\subset \mathbf{L}$ is negative definite and $\sigma_-(\mathbf{L})=19$, one has $\mu(X)\leq 19$ (see \cite{Urabe2}, \emph{cf.}, \cite{Persson}). If $\mu(X)=19$, the quartic is called  \textit{maximizing}.

Given a root lattice $S\subset\mathbf{L}$, let $\tilde{S}_h:=(S_h\otimes\mathbb{Q})\cap\mathbf{L}$ be the primitive hull of $S_h:=S\oplus\Z h$. Since $\sigma_+ \tilde{S}_h^\bot =2$, all positive definite $2$-subspaces in $ \tilde{S}_h^\bot\oplus \mathbb{R}$ can be oriented in a coherent way. Let $\omega$ be one of these coherent orientations. The following statement gives a criterion for the realizability of the triple $( S, h,\mathbf{L} )$ by a simple quartic $X\in \mathbb{P}^3$. It is a combination of the Saint-Donat's description \cite{Donat} of projective models of $K3$-surfaces and the results of Urabe~\cite{Urabe2}.
\begin{proposition} A triple $( S, h,\mathbf{L} )$ is realizable by a simple quartic $X\in \mathbb{P}^3$ (with set of singularities $S$) if and only if
the following conditions satisfied:
\begin{enumerate}
	\item each vector  $e\in (S\otimes\mathbb{Q})\cap \tilde{S}_h$ with $e^2=-2$ and $e\cdot h=0$ lies in $S$,
	\item 	there is no vector $e\in \tilde{S}_h$ such that $e^2=0$ and $e\cdot h=2$
	\end{enumerate}
Then the oriented 2-subspace $\omega_X$ defines the orientation $\omega$. 
\end{proposition}

\subsection{Configurations and $L$-realizations}\label{configuration}
Isomorphism classes of simple singularities are known to be in one-to-one correspondence with those of irreducible root lattices (see Dufree \cite{Dufree}). Hence a set of simple singularities can be identified with a root lattice, the irreducible summands of the latter correspond to the individual singularity points. Thus, the set of simple singularities of a quartic surface $X\subset \mathbb{P}^3$ can be seen as a root lattice $S\subset\mathbf{L}$.

\begin{definition}\label{configuration.admisibility}
	A \textit{configuration} is a finite index extension $\tilde{S}_h\supset S_h=S\oplus\mathbb{Z}h$, $h^2=4$, satisfying the following conditions:
	\begin{enumerate}
		\item each root  $r\in (S\otimes\mathbb{Q})\cap\tilde{S}_h$ with $r^2=-2$ is in $S$,
		\item $\tilde{S}_h$ does not contain an element $v$ with $v^2=0$ and $v\cdot h=2$.
	\end{enumerate}
\end{definition}
An \emph{isomorphism} between two configurations $\tilde{S}'_h,\tilde{S}''_h\supset S_h$ is an isometry  $\tilde{S}'_h \rightarrow \tilde{S}''_h$ preserving both $h$ and $S$ (as a set). We denote by $\Aut_h(\tilde{S}_h)$ the group of \emph{automorphisms} of a configuration $\tilde{S}_h$, \emph{i.e.} autoisometries of $\tilde{S}_h$ preserving $h$. Since $S$ is a characteristic sublattice of  $\tilde S=h^{\perp}_{\tilde{S}_h}$, any isometry of $\tilde S_h$ preserving $h$ preserves $S$; then by item (1) in the Definition \ref{configuration.admisibility}, we have  $\Aut_h(\tilde{S}_h)\subset O(S)$.
\begin{definition}\label{Lrealization}
	An $\mathbf{L}$-\emph{realization} of a configuration $\tilde{S}_h $ is a primitive isometry $\tilde S_h \into \L$.
\end{definition}

Two $\L$-realizations $\tilde{S}'_{h'},\tilde{S}''_{h''} \into \L$ are said to be \emph{isomorphic} if there is an element of the group $O(\L)$ taking $h'$ to $h''$ and $S'$ to $S''$ (as a set). Let $\omega'$ and $\omega''$ be the orientations of these two $\L$-realizations, then these oriented $\L$-realizations are called \emph{strictly isomorphic} if there is an isomorphism between them taking $\omega'$ to $\omega''$, 
An $\L$-realization $\tilde{S}_h \into \L$ is called \emph{symmetric} if it is preserved by an element $a\in O_h(\L)\smallsetminus O_h^+(\L)$, \emph{ i.e.} an autoisometry of $\L$ preserving $S$ (as a set) and $h$ and reversing the positive sign structure; such autoisometries are called as \emph{skew-automorphisms} of the $\L$-realization. If an $\L$-realization  $\tilde{S}_h \into \L$ admits an involutive skew-automorphism, it is called \emph{reflexive}. The notion of isomorphism classes, where we ignore the orientations, may be needed to simplify the classification of $\L$-realizations; namely, we have the following remark.
\begin{remark}
Each isomorphism class consists of one or two strict isomorphism classes depending on whether the $\L$-realizations are symmetric or not, respectively.
\end{remark}
\subsection{Perturbations}\label{perturbations}
Recall that a set of simple singularities can be identified with a root lattice. A \emph{perturbation} of a set of singularities $S$ is a primitive root sublattice $S'$ of $S$. According to E. Looijenga \cite{Looi}, deformation classes of perturbations of an individual simple singular point of type $S$ are in a one-to-one correspondence with the  isomorphism classes of primitive extensions $S'\hookrightarrow S$ of root lattices, see \S \ref{root.lattices} and \S\ref{lattice.extensions}. As shown in \cite{Dynkin}, $S$ admits a perturbation to $S'$
if and only if the Dynkin graph of $S'$ is an induced subgraph of that of $S$. Hence, given a simple quartic $X$, any perturbation of $X$ to a simple quartic $X'$ gives rise to a perturbation of the set of singularities $S$ of $X$ to the set of singularities $S'$ of $X'$. 
According to \cite{Alex.irr.sextics}, the converse also holds: Given a simple quartic
surface $X$ with set of singularities $S$, any perturbation of $S$ to $S'$ is realized by a perturbation of $X$ to $X'$ whose set of singularities is $S'$.

Note that, a perturbation $S'\subset S$ of root lattices give rise to a perturbation of configurations $\tilde{S'}_h\subset\tilde{S}_h$, \ie, a primitive sub-configuration. Here the isotropic subgroup $\K$ is inherited automatically. Thus, any $\L$-realization of $S$ gives a canonical $\L$-realization  of $S'$ as one has the chain of primitive extensions $\tilde{S}'_h\subset\tilde{S}_h\subset\L$.

\subsection{The arithmetical reduction}
Two simple quartics $X_0$ and $X_1$ in  $\mathbb{P}^3$ said to be \emph{equisingular deformation equivalent} if there exists a path $X_t$, $t\in [0.1]$ in the space of simple quartics such that the Milnor number $\mu(X)$ in $X_t$ remains constant. The deformation classification of simple quartics is based on the following statement.


\begin{theorem}[see Theorem 2.3.1 in \cite{AI}]\label{def.class1}
	The map sending a simple quartic surface $X\subset\mathbb{P}^3$ to its oriented $\L$-realization establishes a one to one correspondence between the set of equisingular deformation classes of quartics and that of strict isomorphism classes of oriented $\L$-realizations. Complex  conjugate quartics have isomorphic $\L$-realizations that differ by the orientations.
\end{theorem}

We denote by $\mathcal{X}(S)$ the equisingular deformation class corresponding to $S$ under the bijection given in Theorem \ref{def.class1}. 

\begin{proposition}\label{symmetric.L.realization}
	Cosider an $\L$-realization extending a fixed set of singularities $S$ and let $\mathcal{X}(S)$ be the equisingular deformation class. Then;
	\begin{itemize}
		\item $\mathcal{X}(S)$ in invariant under complex conjugation if and only if $\L$-realization is symmetric
		\item $\mathcal{X}(S)$ contains a real quartic if and only if $\L$-realization is reflexive.
	\end{itemize}
\end{proposition}
According to Proposition \ref{symmetric.L.realization}, symmetric $L$-realizations corresponds to \emph{real}, \emph{i.e.}, conjugation invariant components of $\mathcal{X}(S)$. 

\section{Deformation Classification-Proof of Theorem \ref{th.classification.quartics}}\label{algorithm.for.classification}
Fix a set of singularities $S$ and consider the corresponding $4$-polarized lattice $S_h=S\oplus \Z$, $h^2=4$. Typically the question whether the moduli space $\mathcal{X}:=\mathcal{X}(S)$ is nonempty depends on the polarized lattice $\tilde{S}_h$ only. To assert that $\mathcal{X}\ne \emptyset$ and that a very general member $X\in \mathcal{X}$ has the desired geometric properties we use the results of Nikulin \cite{Niku2} and Saint-Donat \cite{Donat} which reduce the problem to a certain set of conditions given in Definition \ref{configuration.admisibility}. According to Theorem \ref{def.class1} and Definition \ref{configuration.admisibility}, a set of singularities $S$ is realized by a simple quartic surface if and only if  a configuration $\tilde{S}_h$ extending $S_h$ admits a primitive isometry $\tilde{S}_h\into \L$. Hence the general case is splitted into two subcases as finite index extensions $S_h\subset \tilde{S}_h$ as in Definition \ref{configuration.admisibility} and primitive extensions $\tilde{S}_h\into \L$. Theorem \ref{def.class1} states that any configuration $\tilde{S}_h$ give rise to a number of nonempty connected strata $\mathcal{X}(S)$ which are in a bijection with the isomorphism classes of primitive isometries   $\tilde{S}_h \into \mathbf{L}$. A  typical member $X\in \mathcal{X}$ has $NS(X)=X$. The non-generic members $X$ for which $NS(X)\ni h$ fails to be a configuration constitute a countable union of divisors the complement of which is still connected.  Hence the applications of Theorem \ref{def.class1} rely upon the following three question:
\begin{enumerate}
	\item find all configurations $\tilde{S}_h$ (up to isomorphism) extending a given 4-polarized lattice $S_h$;
	\item detect if $\tilde{S}_h$ admits an $\L$-realization;
	\item list all equivalence(isomorphism) classes of $\L$-realizations of  $\tilde{S}_h$.
\end{enumerate}

The classification of $\L$-realizations of configurations $\tilde{S}_h$ extending $S_h$ is done in $4$ steps answering the three questions listed above:

\subsection{Step 1: Enumerating the configurations $\tilde{S}_h$ extending $S_h$}\label{conf}
Question (1) above is settled by Theorem \ref{L-K}, where $\tilde{S}_h$ is determined by a choice of an isotropic subgroup $\K\subset \discr S_h$ where we have $\discr \tilde{S}_h=\K^\bot/\K$. The connected components of the moduli space $\mathcal{X}(S)$ modulo complex conjugation $\operatorname{conj}\colon \mathbb{P}^3\rightarrow\mathbb{P}^3$ are enumerated by the kernel $\K$ of the finite index extension $S_h\subset\tilde{S}_h$ in the given isomorphism class.
\begin{remark}
	If $\K=0$, then $\tilde{S}_h=S_h$ and one has $\discr \tilde{S}_h=\discr S\oplus \langle\frac{1}{4}\rangle$ and $\Aut_h(S_h)=O(S)$. This case corresponds to the classification of nonspecial quartics handled in \cite{Cisem1}.
\end{remark}
There are examples, see Example \ref{multiple.configurations} below, where the set of singularities $S$ admits more than one configuration $\tilde{S}_h$.

\subsection{Step 2 : Detecting if $\tilde{S}_h$ admits an $\L$-realization}\label{existence}
Question (2), which reduces to deciding if genus $g(\tilde{S}_h^{\bot})\neq0$  (in $\mathbf{L} $) in Corollary \ref{bicoset}, is settled by  the existence criterion giving in Theorem \ref{th.N.existence}.
For the first part of the statement, if suffices to list (using Theorem \ref{th.N.existence}) all configurations $ \tilde{S}_h$ that  extends to an $\L$-realization. Implementing the algorithms given in sections \S \ref{conf} and \S \ref{existence} in GAP, we found that there are $4469$ sets of realizable set of singularities, splitting into $8845$ configurations, where $278$ sets of singularities splitting into $347$ configurations are realized by a maximizing quartic. The discussion on perturbation  in the first part the statement is given in \S\ref{Step4} Step 4 below.

\subsection{Step 3 : Listing all isomorphism classes of $\L$-realizations of $\tilde{S}_h$}\label{Step3}
The key to question (3) is the Corollary \ref{bicoset} and to apply it we need tools to list all classes $T\in g(\tilde{S}_h^\bot)$ (assuming the latter is non-empty) and to compute (the cokernels of) the natural homomorphisms 
$$d_S \colon O_h(\tilde{S}_h)\rightarrow \Aut(\operatorname{discr} \tilde{S}_h),\quad d_T \colon O(T)\rightarrow \Aut(\operatorname{discr} T).$$
There are examples, see Example \ref{multiple.Tperb} or Example \ref{multiple.Lrealization3} below, where the genus $g(\tilde{S}_h^\bot)$ does contain more than one isomorphism class, or for a fixed representative $T\in g(\tilde{S}_h^\bot)$,  see Examples \ref{multiple.Lrealization1}, \ref{multiple.Lrealization2} and Example \ref{multiple.Lrealization4} below, where the quotient set given in Corollary \ref{bicoset} does consist of more than one bi-coset, thus giving rise to more than one $\L$-realization. Once the lattice $T=\tilde{S}_h^\bot$ is chosen, one can fix an anti-isometry $\discr \tilde{S}_h \rightarrow \discr T$, and, hence, an isomorphism $\Aut \discr \tilde{S}_h = \Aut \discr T$.

We investigate the isomorphism classes of $\L$-realizations of $\tilde{S}_h$ (\ie, primitive isometries $\tilde{S}_h \into \mathbf{L}$) separately for the maximizing case, \emph{i.e.} $\mu(S)=19$ and non-maximizing case \emph{i.e.} $\mu(S)\le 18$ where for the lattice $T=\tilde{S}_h^\bot$ we use either
\begin{itemize}
	\item Gauss's theory of binary forms \cite{Gauss}, if $T$ is definite of small rank, or
	\item Miranda--Morrison's theory \cite{MM1, MM2, MM3}, if $T$ is indefinite.
\end{itemize}
If $\mu(\mathbf{S})=19$,  the lattice $T=\tilde{S}_h^\bot$ is a positive definite sublattice of rank $2$, and the numbers $(r,c)$ of connected components of the space $\mathcal{X}(S)$ listed in Table \ref{table:maxtable} can easily be computed by Gauss theory of binary quadratic forms \cite{Gauss}. Thus, throughout the rest of the proof we assume $\mu(S)\leq 18$.

If $\mu(\mathbf{S})\le 18$, then $T$ is an indefinite lattice of rank $\operatorname{rk} T \geq 3$, hence  Miranda--Morrison's theory \cite{MM1, MM2, MM3} applies, see \S \ref{MM},  and gives us both $g(\tilde{S}_h^\bot)$ and $\operatorname{Coker}\;d_T$ with in the single finite abelian group $E(T)$ and the natural homomorphism
\begin{align}\label{MMmape}
	 \mathrm{e}\colon \Aut(\operatorname{discr} T)\rightarrow E(T).
\end{align}
Thus, with $\K$, and hence $\tilde{S}_h$ fixed, the further primitive extensions $\tilde{S}_h\into \L$ are enumerated by the cokernel of the well-defined homomorphism
\begin{align}\label{dperb2}
d^\bot\colon\Aut_h(\tilde{S}_h)\rightarrow E(T),
\end{align}
see \S \ref{MM}. In the special case $\K=0$, due to isomorphism $\Aut_h\tilde{S}_h=O(S)$, we have a canonical bijection
$$\pi_0(\mathcal{X}_1(S)/\operatorname{conj})=\operatorname{Coker}[d^\bot\colon O(S)\rightarrow E(T)].$$
assuming that $\tilde{S}_h=S_h$ does admit a primitive extension to $\mathbf{L}$ and taking for $T$ any representative of the genus $S_h^\bot$.

Hence, the $\L$-realization $\tilde{S}_h \into \mathbf{L}$ is unique up to isomorphism, \ie , the space $\mathcal{X}(S)/\operatorname{conj}$ is connected, if and only if the map $d^\bot$ is surjective. 


Distinguishing between a component and its complex conjugate, by Proposition \ref{symmetric.L.realization}, for a fixed $\tilde{S}_h$, the component of the strata $\mathcal{X}(S)$ realizing $\tilde{S}_h$ is  real if and only if the corresponding $\L$-realization extending $\tilde{S}_h$ is symmetric, otherwise it consist of two complex conjugate components  (Asymmetric $\L$-realizations exists, see Example \ref{multiple.configurations}). Thus, to enumerate the real and complex conjugate components of a strata $\mathcal{X}(S)$, one can recast \eqref{MMmape} and \eqref{dperb2} by replacing $\mathrm{e}$ with $\mathrm{e}^+$ and $E(T)$ with $E^+(T)$ and conclude that for each configuration $\tilde{S}_h$,  the corresponding component of the strata $\mathcal{X}(S)$ is real if and only if   
\begin{align}\label{dperb2plus}
d^\bot_+\colon\Aut_h(\tilde{S}_h)\rightarrow E^+(T),
\end{align}
is surjective. It is this latter statement that we proved by computer aided calculations in GAP.

\subsection{Step 4 : Perturbations}\label{Step4}
We compute all iterated perturbations, \ie, perturbations of all families and find out that the $390$ maximizing families given in Table \autoref{table:maxtable} and $39$ extremal families given in Table \autoref{table:extremal39} are not perturbations of anything bigger.

We have effectively implemented all algorithms described in this section in GAP \cite{GAP} and obtain conclusive results that gives us for each realizable set of singularity $S$, the isomorphism classes of configurations $\tilde{S}_h$ extending $S_h$ and for each verified configuration $\tilde{S}_h$,  the strict isomorphism classes of primitive isometries $\tilde{S}_h \into \mathbf{L}$.

In conclusion, the implemented calculations by GAP completes the proof by first listing all the configurations  $\tilde{S}_h$ extending $S_h$ and admitting an $\L$-realization  and then enumerating the strict isomorphism classes of $\L$-realizations, \ie, computing the numbers $(r,c)$ of the stratum  $\mathcal{X}(S)$.

\subsection{Examples}
In this section we illustrate some interesting examples  and demonstrate the calculations for $S=\mathbf{A}_{15}\oplus\mathbf{A}_3$
\begin{example}\label{multiple.configurations}
	The set of singularities $S=2\mathbf{A}_7\oplus\mathbf{A}_3\oplus\mathbf{A}_1$ admits $12$ different configurations $\tilde{S}_h$ (up to isomorphism), Each of them extends to a unique (up to conjugation) $\L$-realization where eleven of them is symmetric and the remaining one is not; hence $\mathcal{X}(S)$ consists of eleven real components and one pair of complex conjugate components
\end{example}

\begin{example}\label{multiple.Tperb}
	The set of singularities $S=\mathbf{A}_{10}\oplus\mathbf{A}_9$ admits a unique configuration $\tilde{S}_h$. It extends to two $\L$-realizations which differ by the lattices $\tilde{S}_h^{\bot}$. One of the $\L$-realizations is symmetric, the other one is not so that $\mathcal{X}(S)$ consists of one real component and one pair of complex conjugate components.
\end{example}

\begin{example}\label{multiple.Lrealization1}
	The set of singularities $S=2\mathbf{D}_{6}\oplus\mathbf{A}_4\oplus\mathbf{A}_3$ admits a unique configuration $\tilde{S}_h$. It extends to two different $\L$-realizations where the two lattices $\tilde{S}_h^{\bot}$ are isomorphic. Both $\L$-realizations are symmetric, hence $\mathcal{X}(S)$ consists of two real components.
\end{example}

\begin{example}\label{multiple.Lrealization2}
	The set of singularities $S=\mathbf{D}_{6}\oplus\mathbf{A}_9\oplus\mathbf{A}_4$ admits a unique configuration $\tilde{S}_h$. It extends to two different $\L$-realizations (with isomorphic lattices $\tilde{S}_h^{\bot}$) where one of them is symmetric and the other one is not; hence $\mathcal{X}(S)$ consists of one real component and one pair of complex conjugate components.
\end{example}

\begin{example}\label{multiple.Lrealization3}
	The set of singularities $S=\mathbf{A}_{9}\oplus\mathbf{A}_6\oplus\mathbf{A}_3\oplus\mathbf{A}_1$ admits two different  configurations $\tilde{S}_h$, each of them extends to two different $\L$-realizations which differ by the lattices $\tilde{S}_h^{\bot}$ where one of them is symmetric and the other one is not; hence $\mathcal{X}(S)$ consists of two real components and two pairs of complex conjugate components..
	\end{example}

\begin{example}\label{multiple.Lrealization4}
The set of singularities $S=\mathbf{A}_{8}\oplus\mathbf{A}_6\oplus\mathbf{A}_3\oplus\mathbf{A}_2$ admits a unique configuration $\tilde{S}_h$. It extends to three different $\L$-realizations (with isomorphic lattices $\tilde{S}_h^{\bot}$) where all of them are not symmetric; hence $\mathcal{X}(S)$ consists of three pairs of complex conjugate components.
\end{example}


\subsection{Demostration for the set of singularites $S=\mathbf{A}_{15}\oplus\mathbf{A}_3$}\label{example.for.proof}

  We demostrate the calculations handled by GAP for the set of singularity $S=\mathbf{A}_{15}\oplus\mathbf{A}_3$. Let $S_h=S\oplus\Z h$, where $h^2=4$.
Then one has
$$\discr S_h\cong\textstyle\langle-\frac{15}{16}\rangle\oplus\langle-\frac{3}{4}\rangle\oplus\langle\frac{1}{4}\rangle \cong (\Z /16\Z)\oplus (\Z /4\Z)\oplus (\Z /4\Z).$$
We fix the generators
$$\text{$\alpha_1$ for $\discr \mathbf{A}_{15}\cong\textstyle\langle-\frac{15}{16}\rangle$,\quad $\alpha_2$ for $\discr \mathbf{A}_3\cong\textstyle\langle-\frac{3}{4}\rangle$,\quad$\alpha_3$ for $\discr \mathbb{Z}h\cong\textstyle\langle\frac{1}{4}\rangle$},$$
and use the coordinate vector notation $[x,y,z]$ for the vector $x\alpha_1+y\alpha_2+z\alpha_3$.

\emph{Step 1}: We determine all isotropic subgroups $\K\subset \discr S_h$ such that the corresponding finite index extension $\tilde{S}_h$ satisfies the conditions in Definition \ref{configuration.admisibility}, \ie, $\tilde{S}_h$ is a configuration extending $S_h$. Up to action of $O(S)$, we have three such isotropic subgroups $\K$ (\ie, three isomorphism classes of configurations $\tilde{S}_h$), which are given in the \autoref{table:kernels}.

\begin{table}
	\centering
	\caption{The isotropic subgroups $\K_i$}\label{table:kernels}
	\begin{tabular}{c c c } 
	\hline \addlinespace[0.11cm]
	 & Generators &  \\ [0.5ex] 
		\hline\addlinespace[0.11cm]
		$\K_1$ & $[8,2,2]$ & cyclic of order $2$  \\ 
		$\K_2$ & $[4,0,2]$ & cyclic of order $4$ \\
	$\K_3$& $[12,2,0]$  & cyclic of order $4$ \\ [1ex] 
		\hline
	\end{tabular}
\end{table}
As the computations given in what follows repeats almost the same for all the configurations, from now on we fix the configuration $\tilde{S}_h$ as the one corresponding to the isotropic subgroup $\K=\K_1$. Then $\discr \tilde{S}_h$, which is given by $\K^\bot/\K$, is generated by the vectors $\{[4,3,1],[4,3,3],[15,3,0]\}$ of orders $2,2$ and $16$, repectively. 

\emph{Step 2}: Consider the orthogonal complement  $T:=\tilde{S}_h^\bot$ given by the  signature $(2, 1)$ and $\discr T =-\discr \tilde{S}_h$. Then by applying Nikulin's Existence Theorem (Theorem \ref{th.N.existence}), we verified that the configuration $\tilde{S}_h$ admits a primitive isometry $\tilde{S}_h \into \mathbf{L}$, \ie, an $\L$-realization.

\emph{Step 3}: The lattice $T=\tilde{S}_h^\bot$ is an indefinite lattice of rank $3$, hence Miranda-Morrison's theory (see \S \ref{MM} and \S\ref{Step3}) can be applied to enumerate the equivalence classes of primitive isometries $\tilde{S}_h \into \mathbf{L}$. By using \eqref{orderofEN}, one gets $|E(T)|=1$ and the map $d^\bot\colon\Aut_h(\tilde{S}_h)\rightarrow E(T)$ is automatically surjective. Thus, as explained in \S \ref{Step3}, the  space $\mathcal{X}(S)/\operatorname{conj}$ is connected.  To enumerate the real and complex conjugate components of $\mathcal{X}(S)$, we need to compute 
$
d^\bot_+\colon\Aut_h(\tilde{S}_h)\rightarrow E^+(T). 
$ (see \eqref{dperbplus}). 
We compute the group $E^+(T)$ directly from the definition given in \eqref{E(N)plus}  which can be restated as
\begin{align}
E^+(T)=\prod_{p|\operatorname{det}(T)}\Gamma_{p,0}/\prod_{p|\operatorname{det}(T)}\Sigma_p^{\sharp}(T)\cdot\Gamma_0^{--}.
\end{align}
Since we have one irregular prime $p=2$, we obtain
$$
E^+(T)=\Gamma_{2,0}/\Sigma_2^{\sharp}(T)\cdot\Gamma_0^{--},
$$
where $\Sigma_2^{\sharp}(T)=\{(1,1),(-1,7),(1,5),(-1,3)\}$. Thus $E^+(T)=\{\pm 1\}$.

The group $\K^\bot/\K$ is given by the following Gram matrix
\[Q:=
\begin{bmatrix} 
1/2 & 0 & 0 \\
0& 1/2 & 0\\
0 & 0 & 5/16 
\end{bmatrix}
\]\label{Q}
in the basis vectors 
$\{\beta_1,\beta_2,\beta_3\}$ where $\beta_1=4\alpha_1+3\alpha_2+\alpha_3$, $\beta_2=4\alpha_1+3\alpha_2+3\alpha_3$ and $\beta_3=15\alpha_1+3\alpha_2$. From now on the coordinate vector notation $[x,y,z]$ will be used for the vector $x\beta_1+y\beta_2+z\beta_3$.  We are interested in the induced action of $\Aut_h(\tilde{S}_h)$ on the discriminant $\discr T=-\K^\bot/\K$. The group $\Aut_h(\tilde{S}_h)$ is generated by a nontrivial symmetry of $\mathbf{A}_{15}$ and a nontrivial symmetry of $\mathbf{A}_3$ where the former give rise to a reflection $t_{a}$ in $\discr T$ with $a=[1,1,2]$ and $a^2=-\frac{9}{4}$. Thus, one has $\mathrm{e^+}(t_a)=1\in E^+(T)$.

A nontrivial symmetry of $\mathbf{A}_3$ induce a reflection $t_{b}$ in $\discr T$ with $b=[1,1,8]$ and $b^2=-21=0 \mod \Z$, \ie, this reflection is one of the exceptional ones listed in  \S\ref{lift}, Proceeding as in \S\ref{lift}, we lift the reflection $t_b$ to $t_{\bar{b}}$ in a $2$-adic lattice $T\otimes\Z_2$:

Since the discriminant $\discr T$ is odd, the algorithm given in \S\ref{lift}, gives us the following two candidate lattices which agrees by Nikulin \cite{Niku2}:
$$
T'=\left[
\begin{array}{rrc}
-2 & 0 & 0 \\
0& -2 & 0\\
0 & 0 & -90
\end{array}
\right]
\quad\mbox{and}\quad
T''=\left[
\begin{array}{rrc}
	-10 & 0 & 0 \\
0& -2 & 0\\
0 & 0 & -90 
\end{array}
\right].
$$
and we choose the second one as the lift since one has $\operatorname{det} T''=-|\discr T|\bmod (\Z_2^{\times})^2$. For the purpose of computing the image of $\mathrm{e}^+$ on the reflection $t_b\in\Aut(\discr T)$ by lifting it to $t_{\bar{b}}$, we replace the first entry of $-Q$ with $-\frac{5}{2}$. Then one gets $\bar{b}^2=-23$, see \eqref{spinn}, and hence,  $\mathrm{e^+}(t_{\bar{b}})=-1\in E^+(T)$. Thus the map $d^\bot_+$ is surjective implying that the strata $\mathcal{X}(S)$ corresponding to $\tilde{S}_h$ consists of one real component.



\subsection{Concluding Remarks}
Clearly any connected component $C\subset \mathcal{X}_*(S)$ containing a real curve is real. However the converse in not true: The known counter examples are discovered in the realm  of irreducible sextics and nonspecial quartics; where the exception for the former is the stratum  $\mathcal{X}_1(\mathbf{A}_7\oplus\mathbf{A}_6\oplus\mathbf{A}_5)$ found by Akyol and Degtyarev \cite{Alex2} and for the latter is the strata $\mathcal{X}_1(\mathbf{A}_7\oplus\mathbf{A}_6\oplus\mathbf{A}_3\oplus\mathbf{A}_2)$ and $\mathcal{X}_1(\mathbf{D}_7\oplus\mathbf{A}_6\oplus\mathbf{A}_3\oplus\mathbf{A}_2)$ found by G\"{u}ne\c{s} Akta\c{s} \cite{aktacs2019real}.
It is worth mentioning that studying phenomena of kind in the whole space of simple sextics or simple quartics would only need an extension of the algorithms, studied in \cite{Alex2} and \cite{aktacs2019real} for both of which one has $\K=0$; the extension should be provided in such a way that  the kernels $\K\neq 0$ are also taken into account. This nontrivial case is still open for simple quartics.

 \subsection{Tables}\label{tables}
This subsection is devoted to present the tables referred in Theorem \ref{th.classification.quartics}.
In \autoref{table:maxtable}, \autoref{table:extremal39} and  \autoref{table:nonreal},
the first column refers the set of simple singularities $S$ which are realized by the families of quartics belonging to the spaces indicated in the table names. For each set of singularity $S$, the column $(r,c)$ gives the numbers of real ($r$) and pairs of complex conjugate ($c$) components of the stratum $\X(S)$ separately for each configuration $\tilde{S}_h$, \ie, for a fixed $S$, the counts $(r,c)$ are aligned in a separate line for each different configuration $\tilde{S}_h$ extending $S_h$. The third column titled as "generators of kernels" gives the description of each separate configuration $\tilde{S}_h$ by listing the generators of the kernel $\K$, as each $\tilde{S}_h$ is determined by a choice of $\K\subset\discr S_h=\discr S\oplus \langle\frac{1}{4}\rangle$. The generators of $\K$ are encoded by using the \emph{glue code} and \emph{glue vectors} introduced in Conway and Sloane \cite{ConvaySloane} (see chapter $16.1$). The only difference it that one more entry for the discriminant group $\langle\frac{1}{4}\rangle$ corresponding to the polarization $h$,  is added to the end of the each glue vector.  To save more space, the brackets of the basis vectors are removed, instead different vectors are separated by a semicolon and only in the cases where a coefficient with two digits appear in the basis vector, comma is used to distinguish the entries. 
\begin{remark}\label{Kummerquartic}
The generators of the kernel of the real equisingular deformation family $\X(S)$ with $S=16\mathbf{A}_1$ listed in the first row of \autoref{table:extremal39} is removed as they are too long to fit in the table. This quartic is known as the \emph{Kummer quartic} described in \cite{Niku.Kummer} which is explicitly generated by 
\begin{align*}
00000111110010002;\,01000011001101110;\,00100101101111000;\\00001101110001110;\, 00011111000110010;\,11110010001110000
\end{align*}
as in the way we display in \autoref{table:extremal39}.
\end{remark}
\begin{remark}\label{3d46a1}
	The generators of the kernel of the real equisingular deformation family $\X(S)$ with $S=3\mathbf{D}_4\oplus6\mathbf{A}_1$ listed in \autoref{table:extremal39} is removed as they are too long to fit in the table. The display of the generators in the  third column is as follows:
	\begin{align*}
0330100012;\, 0110011002;\, 2200001012;\, 3010000112;\, 0001111112
	\end{align*}
\end{remark}
\begin{longtable}{lll}
	\caption{\small The space $\mathcal{X}(S)$ with $\mu(S)=19$} \label{table:maxtable} \\
	\addlinespace[0.11cm]
	\hlineB{2}
	\addlinespace[0.12cm]
	\multicolumn{1}{l}{ Singularities} & \multicolumn{1}{c}{$(r,c)$} &      
	
	\multicolumn{1}{l}{Generators of kernels} \\ 
	\addlinespace[0.12cm]
	\hlineB{2}
	\addlinespace[0.12cm]
	\endfirsthead
	
	\multicolumn{3}{c}%
	{{ \tablename\ \thetable{} -- continued from previous page}} \\
	\addlinespace[0.12cm]
	\hlineB{2} \addlinespace[0.12cm]
	\multicolumn{1}{l}{ Singularities} & \multicolumn{1}{c}{$(r,c)$} & \multicolumn{1}{l}{Generators of kernels} \\ 
	\hlineB{2}
	\addlinespace[0.10cm]
	\endhead
	
	\hline \multicolumn{3}{l}{{Continued on next page}} \\ \hline
	\endfoot
	\addlinespace[0.12cm]
	\hlineB{2.2}
	\addlinespace[0.12cm]
	\endlastfoot	
	\input{max_table}
\end{longtable}

\begin{longtable}{lll}
	\caption{\small Extremal families} \label{table:extremal39} \\
	\addlinespace[0.10cm]
	\hlineB{2.2}
	\addlinespace[0.10cm]
	\multicolumn{1}{l}{ Singularities} & \multicolumn{1}{c}{$(r,c)$} &      
	
	\multicolumn{1}{l}{Generators of kernels} \\ 
	\addlinespace[0.10cm]
	\hlineB{2}
	\addlinespace[0.10cm]
	\endfirsthead
	
	\multicolumn{3}{c}%
	{{\tablename\ \thetable{} -- continued from previous page}} \\
	\addlinespace[0.12cm]
	\hlineB{2} \addlinespace[0.12cm]
	\multicolumn{1}{l}{ Singularities} & \multicolumn{1}{c}{$(r,c)$} & \multicolumn{1}{l}{Generators of kernels} \\ 
	\hlineB{2}
	\addlinespace[0.12cm]
	\endhead
	
	\hline \multicolumn{3}{l}{{Continued on next page}} \\ \hline
	\endfoot
	\addlinespace[0.12cm]
	\hlineB{2}
	\addlinespace[0.12cm]
	\endlastfoot	
	\input{extremal39_table}
	\end{longtable}

\begin{longtable}{lll}
	\caption{\small Nonreal strata  $\mathcal{X}(S)$ with $\mu(S)\le18$} \label{table:nonreal} \\
	\addlinespace[0.12cm]
	\hlineB{2}
	\addlinespace[0.10cm]
	\multicolumn{1}{l}{ Singularities} & \multicolumn{1}{c}{$(r,c)$} &      
	
	\multicolumn{1}{l}{Generators of kernels} \\ 
	\addlinespace[0.12cm]
	\hlineB{2}
	\addlinespace[0.10cm]
	\endfirsthead
	
	\multicolumn{3}{c}%
	{{ \tablename\ \thetable{} -- continued from previous page}} \\
	\addlinespace[0.12cm]
	\hlineB{2} \addlinespace[0.12cm]
	\multicolumn{1}{l}{ Singularities} & \multicolumn{1}{c}{$(r,c)$} & \multicolumn{1}{l}{Generators of kernels} \\ 
	\hlineB{2}
	\addlinespace[0.10cm]
	\endhead
	
	\hline \multicolumn{3}{l}{{Continued on next page}} \\ \hline
	\endfoot
	\addlinespace[0.10cm]
	\hlineB{2}
	\addlinespace[0.10cm]
	\endlastfoot	
	\input{nonreal_table}
\end{longtable}

\bibliographystyle{amsplain}
\bibliography{mybibliography}

\end{document}

%% file: max_table.tex
$\singset{6A3 + A1}$ & $( 1, 0 )$&13030001; 21313000; 30013302\cr
$\singset{A4 + 5A3}$ & $( 1, 0 )$&0312211; 0003313\cr
$\singset{4A4 + A2 + A1}$ & $( 1, 0 )$&3344000\cr
$\singset{A5 + A4 + 3A3 + A1}$ & $( 1, 0 )$&3013311; 3002012\cr
$\singset{3A5 + 4A1}$ & $( 1, 0 )$&30010112; 03310010; 30301010; 44200000\cr
$\singset{3A5 + 2A2}$ & $( 1, 0 )$&303002; 204220; 444000\cr
$\singset{3A5 + A4}$ & $( 1, 0 )$&03302; 24200\cr
$\singset{A6 + A4 + 3A3}$ & $( 1, 0 )$&003331\cr
$\singset{A6 + 2A4 + A3 + A2}$ & $( 2, 0 )$&\cr
$\singset{A6 + A5 + A4 + A3 + A1}$ & $( 2, 0 )$&030212\cr
$\singset{A6 + 2A5 + A3}$ & $( 1, 0 )$&03320\cr & $( 1, 0 )$&03302\cr
$\singset{2A6 + A4 + A2 + A1}$ & $( 0, 1 )$&\cr
$\singset{2A6 + A5 + A2}$ & $( 2, 0 )$&\cr
$\singset{3A6 + A1}$ & $( 2, 0 )$&53600\cr
$\singset{A7 + 3A3 + A2 + A1}$ & $( 1, 0 )$&0331003; 2031002\cr
$\singset{A7 + 4A3}$ & $( 0, 1 )$&602233; 003313\cr & $( 1, 0 )$&222031; 003313\cr & $( 1, 0 )$&031112; 601302\cr
$\singset{A7 + A4 + 2A3 + 2A1}$ & $( 1, 0 )$&2031110; 4000112\cr & $( 1, 0 )$&2001003; 0022112\cr
$\singset{A7 + A4 + 2A3 + A2}$ & $( 1, 0 )$&201102\cr & $( 0, 1 )$&200301\cr
$\singset{A7 + 2A4 + A2 + 2A1}$ & $( 1, 0 )$&4000112\cr
$\singset{A7 + 2A4 + A3 + A1}$ & $( 1, 0 )$&600103\cr
$\singset{A7 + A5 + 2A3 + A1}$ & $( 2, 0 )$&600103; 430010\cr & $( 1, 0 )$&601302; 430010\cr
$\singset{A7 + A5 + A4 + 3A1}$ & $( 1, 0 )$&0301112; 4300010\cr
$\singset{A7 + A5 + A4 + A3}$ & $( 1, 0 )$&20033\cr
$\singset{A7 + A6 + A3 + A2 + A1}$ & $( 2, 0 )$&603001\cr
$\singset{A7 + A6 + 2A3}$ & $( 0, 2 )$&60033\cr & $( 0, 1 )$&60112\cr
$\singset{A7 + A6 + A4 + 2A1}$ & $( 1, 0 )$&400112\cr
$\singset{A7 + A6 + A4 + A2}$ & $( 0, 1 )$&\cr
$\singset{A7 + A6 + A5 + A1}$ & $( 0, 1 )$&40310\cr
$\singset{A7 + 2A6}$ & $( 0, 2 )$&\cr
$\singset{2A7 + A2 + 3A1}$ & $( 1, 0 )$&1501111; 4000112\cr
$\singset{2A7 + A3 + 2A1}$ & $( 0, 1 )$&621103; 042110\cr & $( 1, 0 )$&770103; 042110\cr & $( 1, 0 )$&771102; 400112\cr & $( 1, 0 )$&710103; 641003\cr
$\singset{2A7 + A3 + A2}$ & $( 0, 1 )$&64303\cr & $( 0, 1 )$&60101\cr & $( 1, 0 )$&64303; 44000\cr & $( 1, 0 )$&64303; 66002\cr
$\singset{2A7 + A4 + A1}$ & $( 1, 0 )$&77011\cr
$\singset{2A7 + A5}$ & $( 1, 0 )$&7333\cr
$\singset{A8 + 3A3 + A2}$ & $( 1, 0 )$&033303\cr
$\singset{A8 + A5 + A3 + A2 + A1}$ & $( 0, 1 )$&032012; 640100\cr
$\singset{A8 + A5 + A4 + A2}$ & $( 2, 0 )$&34020\cr
$\singset{A8 + A6 + A3 + A2}$ & $( 0, 3 )$&\cr
$\singset{A8 + A6 + A4 + A1}$ & $( 0, 1 )$&\cr
$\singset{A8 + A6 + A5}$ & $( 1, 1 )$&\cr
$\singset{A8 + A7 + A3 + A1}$ & $( 2, 0 )$&06303\cr
$\singset{2A8 + A2 + A1}$ & $( 1, 1 )$&33000\cr
$\singset{A9 + A4 + A3 + A2 + A1}$ & $( 1, 1 )$&502010\cr & $( 1, 1 )$&500012\cr
$\singset{A9 + 2A4 + 2A1}$ & $( 1, 0 )$&500012; 644000\cr
$\singset{A9 + 2A4 + A2}$ & $( 1, 0 )$&41100\cr
$\singset{A9 + A5 + A3 + A2}$ & $( 2, 0 )$&53000\cr
$\singset{A9 + 2A5}$ & $( 0, 1 )$&5030\cr & $( 1, 0 )$&0332\cr
$\singset{A9 + A6 + A2 + 2A1}$ & $( 1, 0 )$&500102\cr
$\singset{A9 + A6 + 2A2}$ & $( 1, 0 )$&\cr
$\singset{A9 + A6 + A3 + A1}$ & $( 1, 1 )$&50210\cr & $( 1, 1 )$&50012\cr
$\singset{A9 + A7 + 3A1}$ & $( 1, 0 )$&501110; 040112\cr
$\singset{A9 + A7 + A2 + A1}$ & $( 0, 1 )$&50012\cr
$\singset{A9 + A7 + A3}$ & $( 2, 0 )$&0633\cr
$\singset{A9 + A8 + 2A1}$ & $( 0, 1 )$&50102\cr
$\singset{A9 + A8 + A2}$ & $( 1, 1 )$&\cr
$\singset{2A9 + A1}$ & $( 1, 0 )$&5012; 8600\cr & $( 1, 0 )$&5502; 8600\cr
$\singset{A10 + 3A3}$ & $( 1, 0 )$&03331\cr
$\singset{A10 + A4 + A3 + A2}$ & $( 0, 1 )$&\cr
$\singset{A10 + A5 + A3 + A1}$ & $( 0, 1 )$&03212\cr
$\singset{A10 + A5 + A4}$ & $( 1, 0 )$&\cr
$\singset{A10 + A6 + A2 + A1}$ & $( 1, 0 )$&\cr
$\singset{A10 + A6 + A3}$ & $( 0, 2 )$&\cr
$\singset{A10 + A7 + 2A1}$ & $( 1, 0 )$&04112\cr
$\singset{A10 + A7 + A2}$ & $( 0, 2 )$&\cr
$\singset{A10 + A8 + A1}$ & $( 0, 1 )$&\cr
$\singset{A10 + A9}$ & $( 1, 1 )$&\cr
$\singset{A11 + A3 + 2A2 + A1}$ & $( 1, 0 )$&300001; 402100\cr & $( 1, 0 )$&310002; 402100\cr
$\singset{A11 + 2A3 + 2A1}$ & $( 1, 0 )$&311011; 022112\cr
$\singset{A11 + 2A3 + A2}$ & $( 1, 0 )$&92201\cr & $( 1, 0 )$&30003\cr & $( 0, 1 )$&33002\cr & $( 0, 1 )$&33200\cr
$\singset{A11 + A4 + 2A2}$ & $( 1, 0 )$&80120\cr & $( 1, 0 )$&60002; 80120\cr & $( 1, 0 )$&30001; 80120\cr
$\singset{A11 + A4 + A3 + A1}$ & $( 1, 0 )$&90001\cr & $( 1, 0 )$&90102\cr
$\singset{A11 + A5 + 3A1}$ & $( 1, 0 )$&031112; 600110; 820000\cr
$\singset{A11 + A5 + A2 + A1}$ & $( 1, 0 )$&93011; 82000\cr & $( 1, 0 )$&90001; 82000\cr
$\singset{A11 + A5 + A3}$ & $( 1, 0 )$&9001\cr & $( 1, 0 )$&9001; 8200\cr & $( 1, 0 )$&3032\cr & $( 1, 0 )$&3032; 8200\cr
$\singset{A11 + A6 + 2A1}$ & $( 0, 1 )$&60110\cr & $( 1, 0 )$&30001\cr
$\singset{A11 + A6 + A2}$ & $( 0, 2 )$&\cr & $( 0, 1 )$&6002\cr & $( 2, 0 )$&9003\cr
$\singset{A11 + A7 + A1}$ & $( 1, 0 )$&9403\cr & $( 1, 0 )$&9003\cr & $( 0, 1 )$&9613\cr
$\singset{A11 + A8}$ & $( 1, 0 )$&903\cr
$\singset{A12 + A3 + 2A2}$ & $( 2, 0 )$&\cr
$\singset{A12 + A4 + A2 + A1}$ & $( 0, 1 )$&\cr
$\singset{A12 + A5 + A2}$ & $( 1, 1 )$&\cr
$\singset{A12 + A6 + A1}$ & $( 1, 1 )$&\cr
$\singset{A13 + A3 + A2 + A1}$ & $( 0, 1 )$&70010\cr
$\singset{A13 + A4 + 2A1}$ & $( 1, 0 )$&70100\cr
$\singset{A13 + A4 + A2}$ & $( 1, 0 )$&\cr
$\singset{A13 + A5 + A1}$ & $( 1, 0 )$&7302\cr & $( 1, 0 )$&7010\cr
$\singset{A13 + A6}$ & $( 0, 2 )$&\cr
$\singset{A14 + 2A2 + A1}$ & $( 1, 0 )$&52000\cr
$\singset{A14 + A3 + A2}$ & $( 0, 2 )$&\cr & $( 0, 1 )$&5020\cr
$\singset{A14 + A5}$ & $( 0, 2 )$&\cr
$\singset{A15 + A2 + 2A1}$ & $( 1, 0 )$&20101\cr
$\singset{A15 + 2A2}$ & $( 0, 1 )$&\cr & $( 1, 0 )$&8000\cr & $( 1, 0 )$&4002\cr
$\singset{A15 + A3 + A1}$ & $( 0, 1 )$&12,3,1,3\cr & $( 1, 0 )$&10,3,1,2\cr & $( 1, 0 )$&10,0,1,3\cr
$\singset{A15 + A4}$ & $( 1, 0 )$&12,0,2\cr
$\singset{A16 + A2 + A1}$ & $( 1, 0 )$&\cr
$\singset{A17 + 2A1}$ & $( 1, 0 )$&9102\cr & $( 1, 0 )$&9102; 6000\cr
$\singset{A17 + A2}$ & $( 1, 1 )$&\cr & $( 1, 0 )$&12,0,0\cr
$\singset{A18 + A1}$ & $( 1, 1 )$&\cr
$\singset{A19}$ & $( 1, 0 )$&10,2\cr
$\singset{D4 + 5A3}$ & $( 1, 0 )$&2102113; 2031110; 2002022\cr
$\singset{D4 + 2A5 + 2A2 + A1}$ & $( 1, 0 )$&2330000; 3300012; 0442200\cr
$\singset{D4 + 2A5 + A3 + 2A1}$ & $( 1, 0 )$&1302010; 2030012; 2300102\cr
$\singset{D4 + 2A5 + A4 + A1}$ & $( 1, 0 )$&333000; 203012\cr
$\singset{D4 + 3A5}$ & $( 1, 0 )$&13030; 23300; 02420\cr
$\singset{D4 + A7 + 2A3 + A2}$ & $( 1, 0 )$&063003; 240200\cr & $( 1, 0 )$&223300; 202202\cr
$\singset{D4 + A7 + A4 + 2A2}$ & $( 1, 0 )$&240002\cr
$\singset{D4 + A7 + A5 + A2 + A1}$ & $( 1, 0 )$&103012; 140002\cr
$\singset{D4 + A7 + A6 + A2}$ & $( 1, 0 )$&24002\cr
$\singset{D4 + 2A7 + A1}$ & $( 1, 0 )$&37511; 34002\cr
$\singset{D4 + A9 + A5 + A1}$ & $( 1, 0 )$&35302; 20312\cr
$\singset{D4 + A11 + 2A2}$ & $( 1, 0 )$&26000; 04210\cr & $( 1, 0 )$&09001; 04210\cr
$\singset{2D4 + 3A3 + A2}$ & $( 1, 0 )$&1313303; 3220200; 1300202\cr
$\singset{4D4 + 3A1}$ & $( 1, 0 )$&21201010;30031102;00110112;23021100;22000112\cr
$\singset{D5 + 2A7}$ & $( 1, 0 )$&1730; 2042\cr
$\singset{D6 + A5 + A4 + A3 + A1}$ & $( 1, 0 )$&300212; 230012\cr & $( 1, 0 )$&300212; 230210\cr
$\singset{D6 + A5 + 2A4}$ & $( 1, 0 )$&33002\cr
$\singset{D6 + 2A5 + A3}$ & $( 1, 0 )$&33002; 23300\cr & $( 1, 0 )$&33020; 23300\cr
$\singset{D6 + A6 + A5 + A2}$ & $( 2, 0 )$&10302\cr
$\singset{D6 + A7 + A3 + A2 + A1}$ & $( 1, 0 )$&063001; 102012\cr
$\singset{D6 + A7 + 2A3}$ & $( 1, 0 )$&22233; 24200\cr & $( 1, 0 )$&06332; 20222\cr
$\singset{D6 + A7 + A4 + A2}$ & $( 1, 0 )$&24002\cr
$\singset{D6 + A7 + A5 + A1}$ & $( 0, 1 )$&34010; 14302\cr & $( 1, 0 )$&24002; 04310\cr
$\singset{D6 + A7 + A6}$ & $( 0, 1 )$&2402\cr
$\singset{D6 + A8 + A5}$ & $( 1, 1 )$&1032\cr
$\singset{D6 + A9 + 2A2}$ & $( 1, 0 )$&15000\cr
$\singset{D6 + A9 + A3 + A1}$ & $( 0, 1 )$&15202; 30212\cr
$\singset{D6 + A9 + A4}$ & $( 1, 1 )$&1500\cr
$\singset{D6 + A11 + A2}$ & $( 0, 1 )$&2600\cr & $( 1, 0 )$&0301\cr
$\singset{D6 + A13}$ & $( 0, 1 )$&172\cr
$\singset{D6 + D4 + A5 + A3 + A1}$ & $( 1, 0 )$&313202; 330012; 110210\cr & $( 1, 0 )$&313000; 330012; 110210\cr
$\singset{D6 + 2D4 + A3 + 2A1}$ & $( 1, 0 )$&1102100; 2120110; 2330002; 2212000\cr
$\singset{D6 + D5 + A7 + A1}$ & $( 1, 0 )$&23213; 12012\cr
$\singset{2D6 + A4 + A3}$ & $( 2, 0 )$&33020; 22022\cr
$\singset{2D6 + A5 + 2A1}$ & $( 1, 0 )$&110110; 033110; 203102\cr
$\singset{2D6 + A5 + A2}$ & $( 1, 0 )$&10302; 33002\cr
$\singset{2D6 + A7}$ & $( 1, 0 )$&1102; 3342\cr
$\singset{2D6 + D4 + 3A1}$ & $( 1, 0 )$&1121102; 2031102; 0310102; 1010012\cr
$\singset{2D6 + D4 + A2 + A1}$ & $( 1, 0 )$&121010; 303012; 031012\cr
$\singset{2D6 + D5 + 2A1}$ & $( 1, 0 )$&212100; 032012; 302102\cr
$\singset{3D6 + A1}$ & $( 1, 0 )$&31200; 23300; 20112\cr & $( 1, 0 )$&31200; 13002; 30302\cr
$\singset{D7 + 4A3}$ & $( 1, 0 )$&110321; 103213\cr & $( 1, 0 )$&113122; 010313\cr
$\singset{D7 + 2A4 + 2A2}$ & $( 1, 0 )$&\cr
$\singset{D7 + A5 + A4 + A2 + A1}$ & $( 0, 1 )$&230012\cr
$\singset{D7 + 2A5 + 2A1}$ & $( 1, 0 )$&033110; 203012\cr
$\singset{D7 + A6 + A4 + A2}$ & $( 0, 1 )$&\cr
$\singset{D7 + A6 + A5 + A1}$ & $( 2, 0 )$&20312\cr
$\singset{D7 + 2A6}$ & $( 0, 1 )$&\cr
$\singset{D7 + A7 + A3 + 2A1}$ & $( 1, 0 )$&123112; 202112\cr & $( 1, 0 )$&120113; 202112\cr & $( 1, 0 )$&063001; 202112\cr
$\singset{D7 + A7 + A3 + A2}$ & $( 0, 1 )$&16100\cr & $( 0, 1 )$&16203\cr & $( 0, 1 )$&02301\cr
$\singset{D7 + A9 + A2 + A1}$ & $( 1, 0 )$&25010\cr & $( 1, 0 )$&05012\cr
$\singset{D7 + A10 + A2}$ & $( 0, 1 )$&\cr
$\singset{D7 + A11 + A1}$ & $( 1, 0 )$&1900\cr & $( 1, 0 )$&0903\cr
$\singset{D7 + D6 + A5 + A1}$ & $( 1, 0 )$&23300; 21012\cr & $( 1, 0 )$&03302; 21012\cr
$\singset{D7 + 2D6}$ & $( 1, 0 )$&2222; 0132\cr
$\singset{2D7 + A3 + A2}$ & $( 1, 0 )$&11101\cr
$\singset{D8 + 3A3 + A2}$ & $( 1, 0 )$&313103; 302002\cr
$\singset{D8 + A5 + A3 + 3A1}$ & $( 1, 0 )$&3301000; 2021012; 1001102\cr & $( 1, 0 )$&3321002; 2021012; 1001102\cr
$\singset{D8 + A5 + A4 + 2A1}$ & $( 1, 0 )$&330100; 100112\cr
$\singset{D8 + A6 + A3 + A2}$ & $( 0, 1 )$&30202\cr
$\singset{D8 + A7 + A2 + 2A1}$ & $( 1, 0 )$&300112; 040112\cr
$\singset{D8 + A7 + A3 + A1}$ & $( 1, 0 )$&36103; 30202\cr
$\singset{D8 + A9 + 2A1}$ & $( 1, 0 )$&10112; 35102\cr
$\singset{D8 + D4 + A3 + 2A2}$ & $( 1, 0 )$&122000; 232002\cr
$\singset{D8 + D4 + A5 + 2A1}$ & $( 1, 0 )$&120110; 310002; 223100\cr
$\singset{D8 + D5 + A5 + A1}$ & $( 1, 0 )$&02312; 30310\cr
$\singset{D8 + D6 + A3 + 2A1}$ & $( 1, 0 )$&110010; 232010; 300112\cr & $( 1, 0 )$&110010; 230012; 302110\cr
$\singset{D8 + D6 + A4 + A1}$ & $( 1, 0 )$&21012; 13010\cr
$\singset{D8 + D6 + A5}$ & $( 1, 0 )$&2130; 1332\cr
$\singset{D8 + D6 + D4 + A1}$ & $( 1, 0 )$&11212; 03112; 30102\cr
$\singset{D8 + D7 + 2A2}$ & $( 1, 0 )$&32002\cr
$\singset{2D8 + 3A1}$ & $( 1, 0 )$&221012; 300112; 010112\cr
$\singset{2D8 + A2 + A1}$ & $( 1, 0 )$&31000; 12002\cr
$\singset{D9 + A5 + A4 + A1}$ & $( 1, 0 )$&23012\cr
$\singset{D9 + 2A5}$ & $( 1, 0 )$&2330\cr
$\singset{D9 + A6 + 2A2}$ & $( 1, 0 )$&\cr
$\singset{D9 + A7 + A2 + A1}$ & $( 1, 0 )$&32013\cr
$\singset{D9 + A9 + A1}$ & $( 0, 1 )$&2510\cr
$\singset{D10 + A4 + 2A2 + A1}$ & $( 1, 0 )$&300012\cr
$\singset{D10 + A4 + A3 + 2A1}$ & $( 2, 0 )$&202112; 300102\cr
$\singset{D10 + A5 + A2 + 2A1}$ & $( 1, 0 )$&130000; 300102\cr
$\singset{D10 + A5 + A3 + A1}$ & $( 1, 0 )$&13202; 30210\cr & $( 1, 0 )$&13000; 30210\cr & $( 1, 0 )$&13202; 30012\cr & $( 1, 0 )$&13000; 30012\cr
$\singset{D10 + A5 + A4}$ & $( 1, 0 )$&3300\cr
$\singset{D10 + A6 + A2 + A1}$ & $( 1, 0 )$&30012\cr
$\singset{D10 + A7 + 2A1}$ & $( 0, 1 )$&24110; 14012\cr
$\singset{D10 + A7 + A2}$ & $( 1, 0 )$&2402\cr
$\singset{D10 + A8 + A1}$ & $( 1, 0 )$&1012\cr
$\singset{D10 + A9}$ & $( 1, 0 )$&152\cr
$\singset{D10 + D6 + A2 + A1}$ & $( 1, 0 )$&33000; 10012\cr
$\singset{D10 + D6 + A3}$ & $( 1, 0 )$&3100; 1322\cr
$\singset{D10 + D7 + 2A1}$ & $( 1, 0 )$&22112; 10102\cr
$\singset{D10 + D8 + A1}$ & $( 1, 0 )$&2302; 1112\cr
$\singset{D11 + 2A3 + A2}$ & $( 1, 0 )$&11101\cr
$\singset{D11 + A6 + A2}$ & $( 0, 1 )$&\cr
$\singset{D11 + A7 + A1}$ & $( 1, 0 )$&3203\cr
$\singset{D12 + A3 + 2A2}$ & $( 1, 0 )$&32000\cr & $( 1, 0 )$&30002\cr
$\singset{D12 + A4 + A2 + A1}$ & $( 1, 0 )$&10002\cr
$\singset{D12 + A5 + 2A1}$ & $( 1, 0 )$&10110; 33102\cr
$\singset{D12 + A5 + A2}$ & $( 1, 0 )$&3002\cr
$\singset{D12 + A6 + A1}$ & $( 0, 1 )$&3002\cr
$\singset{D12 + D6 + A1}$ & $( 1, 0 )$&1200; 3112\cr
$\singset{D13 + A5 + A1}$ & $( 1, 0 )$&2312\cr
$\singset{D14 + A3 + 2A1}$ & $( 1, 0 )$&22112; 12012\cr
$\singset{D14 + A4 + A1}$ & $( 1, 0 )$&1010\cr
$\singset{D14 + A5}$ & $( 1, 0 )$&132\cr
$\singset{D15 + 2A2}$ & $( 1, 0 )$&\cr
$\singset{D16 + A2 + A1}$ & $( 1, 0 )$&3000\cr
$\singset{D18 + A1}$ & $( 1, 0 )$&312\cr
$\singset{E6 + A11 + 2A1}$ & $( 1, 0 )$&03001; 18000\cr
$\singset{E6 + A11 + A2}$ & $( 1, 0 )$&0301; 1800\cr
$\singset{E6 + A12 + A1}$ & $( 1, 0 )$&\cr
$\singset{E6 + A13}$ & $( 1, 0 )$&\cr
$\singset{E6 + D8 + D5}$ & $( 1, 0 )$&0122\cr
$\singset{E6 + D9 + A4}$ & $( 1, 0 )$&\cr
$\singset{E6 + D10 + A2 + A1}$ & $( 1, 0 )$&01012\cr
$\singset{E6 + D12 + A1}$ & $( 1, 0 )$&0302\cr
$\singset{E6 + D13}$ & $( 1, 0 )$&\cr
$\singset{2E6 + D7}$ & $( 1, 0 )$&\cr
$\singset{3E6 + A1}$ & $( 1, 0 )$&21100\cr
$\singset{E7 + 2A4 + A3 + A1}$ & $( 1, 0 )$&100212\cr
$\singset{E7 + A5 + A4 + A3}$ & $( 1, 0 )$&13002\cr & $( 1, 0 )$&13020\cr
$\singset{E7 + A6 + A4 + A2}$ & $( 1, 0 )$&\cr
$\singset{E7 + A6 + A5 + A1}$ & $( 1, 0 )$&10302\cr
$\singset{E7 + 2A6}$ & $( 0, 1 )$&\cr
$\singset{E7 + A7 + A3 + 2A1}$ & $( 1, 0 )$&063001; 102102\cr
$\singset{E7 + A7 + A3 + A2}$ & $( 1, 0 )$&06103\cr
$\singset{E7 + A7 + A4 + A1}$ & $( 0, 1 )$&14010\cr
$\singset{E7 + A7 + A5}$ & $( 0, 1 )$&1032\cr
$\singset{E7 + A8 + A3 + A1}$ & $( 0, 1 )$&10212\cr
$\singset{E7 + A8 + A4}$ & $( 2, 0 )$&\cr
$\singset{E7 + A9 + A2 + A1}$ & $( 1, 0 )$&15000\cr & $( 1, 0 )$&05012\cr
$\singset{E7 + A9 + A3}$ & $( 0, 1 )$&1500\cr
$\singset{E7 + A10 + A2}$ & $( 0, 1 )$&\cr
$\singset{E7 + A11 + A1}$ & $( 1, 0 )$&1911\cr & $( 1, 0 )$&0903\cr
$\singset{E7 + A12}$ & $( 1, 1 )$&\cr
$\singset{E7 + D4 + A5 + A2 + A1}$ & $( 1, 0 )$&113000; 130012\cr
$\singset{E7 + D4 + A7 + A1}$ & $( 1, 0 )$&10410; 12012\cr
$\singset{E7 + D5 + A5 + 2A1}$ & $( 1, 0 )$&023012; 120102\cr
$\singset{E7 + D5 + A7}$ & $( 1, 0 )$&1323\cr
$\singset{E7 + D6 + A4 + A2}$ & $( 1, 0 )$&13002\cr
$\singset{E7 + D6 + A5 + A1}$ & $( 1, 0 )$&03302; 11002\cr
$\singset{E7 + D6 + A6}$ & $( 0, 1 )$&1102\cr
$\singset{E7 + D6 + D5 + A1}$ & $( 1, 0 )$&12012; 01212\cr
$\singset{E7 + 2D6}$ & $( 1, 0 )$&0132; 1120\cr
$\singset{E7 + D7 + A4 + A1}$ & $( 2, 0 )$&12012\cr
$\singset{E7 + D7 + A5}$ & $( 1, 0 )$&1230\cr & $( 1, 0 )$&1032\cr
$\singset{E7 + D8 + A2 + 2A1}$ & $( 1, 0 )$&030112; 120102\cr
$\singset{E7 + D8 + A3 + A1}$ & $( 1, 0 )$&01202; 11010\cr
$\singset{E7 + D10 + 2A1}$ & $( 1, 0 )$&11000; 03012\cr
$\singset{E7 + D10 + A2}$ & $( 1, 0 )$&1100\cr
$\singset{E7 + D11 + A1}$ & $( 1, 0 )$&1212\cr
$\singset{E7 + D12}$ & $( 1, 0 )$&012\cr
$\singset{E7 + E6 + A6}$ & $( 1, 0 )$&\cr
$\singset{E7 + E6 + D5 + A1}$ & $( 1, 0 )$&10212\cr
$\singset{E7 + E6 + D6}$ & $( 1, 0 )$&1012\cr
$\singset{2E7 + A4 + A1}$ & $( 1, 0 )$&11002\cr
$\singset{2E7 + A5}$ & $( 1, 0 )$&1032\cr & $( 1, 0 )$&1102\cr
$\singset{2E7 + D4 + A1}$ & $( 1, 0 )$&11300; 10112\cr
$\singset{2E7 + D5}$ & $( 1, 0 )$&1102\cr
$\singset{E8 + 3A3 + A2}$ & $( 1, 0 )$&031101\cr
$\singset{E8 + 2A4 + A2 + A1}$ & $( 1, 0 )$&\cr
$\singset{E8 + A6 + A3 + A2}$ & $( 0, 1 )$&\cr
$\singset{E8 + A6 + A4 + A1}$ & $( 1, 0 )$&\cr
$\singset{E8 + A6 + A5}$ & $( 1, 0 )$&\cr
$\singset{E8 + A7 + A2 + 2A1}$ & $( 1, 0 )$&040112\cr
$\singset{E8 + A7 + A3 + A1}$ & $( 1, 0 )$&02303\cr
$\singset{E8 + A9 + 2A1}$ & $( 1, 0 )$&05102\cr
$\singset{E8 + A9 + A2}$ & $( 1, 0 )$&\cr
$\singset{E8 + A10 + A1}$ & $( 1, 0 )$&\cr
$\singset{E8 + A11}$ & $( 1, 0 )$&091\cr
$\singset{E8 + D5 + A5 + A1}$ & $( 1, 0 )$&02312\cr
$\singset{E8 + D6 + A5}$ & $( 1, 0 )$&0132\cr
$\singset{E8 + D7 + 2A2}$ & $( 1, 0 )$&\cr
$\singset{E8 + D10 + A1}$ & $( 1, 0 )$&0312\cr
$\singset{E8 + E6 + A4 + A1}$ & $( 1, 0 )$&\cr
$\singset{E8 + E6 + D5}$ & $( 1, 0 )$&\cr
$\singset{E8 + E7 + A3 + A1}$ & $( 1, 0 )$&01212\cr
$\singset{E8 + E7 + A4}$ & $( 1, 0 )$&\cr
$\singset{2E8 + A2 + A1}$ & $( 1, 0 )$&

%% file: extremal39_table.tex
$\singset{16A1}$ & $( 1, 0 )$&see Remark \ref{Kummerquartic}\cr
$\singset{4D4}$ & $( 1, 0 )$&13320; 22110\cr
$\singset{D4 + 3A3 + 4A1}$ & $( 1, 0 )$&133101103; 222011000; 002211110\cr
$\singset{6A3}$ & $( 1, 0 )$&1320330; 2031130\cr & $( 1, 0 )$&1313002; 3003332; 0022022\cr
$\singset{2A7 + 4A1}$ & $( 1, 0 )$&2600110; 0411110\cr
$\singset{2D4 + 2A3 + 4A1}$ & $( 1, 0 )$&202010102; 310211000; 222200000; 002211110\cr
$\singset{3D4 + 6A1}$ & $( 1, 0 )$&see Remark \ref{3d46a1}\cr
$\singset{D5 + 3A3 + 4A1}$ & $( 1, 0 )$&231100113; 202211000; 202000112\cr
$\singset{D5 + D4 + 2A3 + 3A1}$ & $( 1, 0 )$&30111113; 03201012; 01220110\cr
$\singset{2D5 + 3A2 + 2A1}$ & $( 1, 0 )$&22000112\cr
$\singset{2D5 + 2A3 + 2A1}$ & $( 1, 0 )$&1331000; 2020112\cr & $( 1, 0 )$&1303113; 2020112\cr
$\singset{2D5 + 2A4}$ & $( 1, 0 )$&\cr
$\singset{2D5 + A7 + A1}$ & $( 1, 0 )$&31200\cr
$\singset{2D5 + A8}$ & $( 1, 0 )$&\cr
$\singset{2D5 + D4 + A4}$ & $( 1, 0 )$&22102\cr
$\singset{2D5 + 2D4}$ & $( 1, 0 )$&02322; 20112\cr
$\singset{3D5 + A2 + A1}$ & $( 1, 0 )$&131013\cr
$\singset{3D5 + A3}$ & $( 1, 0 )$&33213\cr
$\singset{D6 + 3A3 + 3A1}$ & $( 1, 0 )$&11330103; 10220100; 20201102\cr
$\singset{D6 + A7 + A3 + 2A1}$ & $( 1, 0 )$&323101; 202112\cr
$\singset{D6 + D4 + 2A3 + 2A1}$ & $( 1, 0 )$&3002012; 1020102; 2200112\cr
$\singset{D6 + D5 + A3 + 4A1}$ & $( 1, 0 )$&32210000; 12001110; 02200112\cr
$\singset{D6 + D5 + 2A3 + A1}$ & $( 1, 0 )$&213113; 102012\cr
$\singset{2D6 + 2A3}$ & $( 1, 0 )$&22220; 33020\cr
$\singset{2D6 + D4 + 2A1}$ & $( 1, 0 )$&312112; 033012; 101012\cr
$\singset{D7 + 2D5 + A1}$ & $( 1, 0 )$&33103\cr
$\singset{D8 + 2D5}$ & $( 1, 0 )$&2222\cr & $( 1, 0 )$&3220\cr
$\singset{D9 + D5 + A3 + A1}$ & $( 1, 0 )$&31303\cr
$\singset{D9 + D5 + D4}$ & $( 1, 0 )$&2222\cr
$\singset{2D9}$ & $( 1, 0 )$&\cr
$\singset{E6 + A5 + 2A2 + 3A1}$ & $( 1, 0 )$&03001112; 14120000\cr
$\singset{E6 + A5 + A3 + 2A2}$ & $( 1, 0 )$&240110\cr
$\singset{E6 + 2A5 + A2}$ & $( 1, 0 )$&22400\cr
$\singset{E6 + A11 + A1}$ & $( 1, 0 )$&1800\cr & $( 1, 0 )$&0602; 1800\cr
$\singset{E6 + D4 + 4A2}$ & $( 1, 0 )$&2022220

%% file: nonreal_table.tex
$\singset{5A3 + A1}$ & $( 0, 1 )$&1331313\cr
$\singset{A7 + 3A3 + A1}$ & $( 0, 1 )$&213113\cr
$\singset{2A6 + 2A3}$ & $( 0, 1 )$&\cr
$\singset{3A6}$ & $( 0, 1 )$&\cr
$\singset{2A7 + 2A2}$ & $( 0, 1 )$&\cr
$\singset{2A7 + A3 + A1}$ & $( 0, 1 )$&22311\cr
$\singset{A11 + A3 + A2 + 2A1}$ & $( 0, 1 )$&600110\cr
$\singset{A11 + 2A3 + A1}$ & $( 0, 1 )$&31111\cr
$\singset{D5 + 2A6 + A1}$ & $( 0, 1 )$&\cr
$\singset{D5 + A9 + A3 + A1}$ & $( 0, 1 )$&25010\cr
$\singset{D6 + 2A6}$ & $( 0, 1 )$&\cr
$\singset{D6 + A7 + A3 + 2A1}$ & $( 0, 1 )$&342102; 102012\cr
$\singset{E7 + A7 + A3 + A1}$ & $( 0, 1 )$&14010